%% file: ex_article.tex
\documentclass[final,hidelinks,onefignum,onetabnum]{siamart251216}

\input{ex_shared}

\ifpdf
\hypersetup{
  pdftitle={Scalable Fully Implicit Runge--Kutta Solvers for
Time-Dependent Rayleigh--B\'enard Convection},
  pdfauthor={Ahsan Ali and Robert C. Kirby}
}
\fi

\usepackage{tikz}
\usetikzlibrary{calc}
\usepackage{booktabs,graphicx,listings}\usepackage[section]{placeins}

\providecommand{\Ra}{\mathrm{Ra}}
\providecommand{\Pran}{\mathrm{Pr}}

\begin{document}

\maketitle

\begin{abstract}
We consider the efficient numerical solution of time-dependent Rayleigh--B\'enard convection governed by the incompressible Boussinesq equations. The coupled velocity--pressure--temperature system presents significant computational challenges arising from incompressibility, buoyancy coupling, convection, and the simultaneous resolution of multiple physical time scales. We present a systematic study of fully implicit
Runge--Kutta methods for time-dependent
Rayleigh--B\'enard convection, together with
scalable finite element solvers for the resulting
stage-coupled systems. We employ high-order, L-stable Radau IIA methods for temporal integration and investigate Taylor--Hood and Alfeld--Sorokina finite elements in two dimensions, along with Taylor--Hood and enriched macroelement discretizations in three dimensions. The macroelement approaches offer strong local inf--sup stability and patch structures well suited to Vanka-type relaxation. 
In two dimensions, monolithic multigrid achieves
nearly stage-independent outer Krylov iteration counts
and maintains robust convergence under mesh refinement. On the same mesh with incremental pressure correction, selected Radau configurations reduce elapsed solve time by 18--32\% and yield smaller sampled transient differences relative to a numerical reference solution. Three-dimensional experiments demonstrate favorable parallel scalability with Taylor--Hood elements and substantially improved discrete mass conservation with the enriched macroelement discretization. These results identify regimes in which coupled multistage fully implicit integration is computationally competitive when equipped with appropriate multigrid preconditioners.
\end{abstract}



\section{Introduction}
\label{sec:introduction}

Rayleigh--B\'enard convection couples incompressible fluid motion to heat transport through buoyancy. Even on a simple domain, the numerical method must accommodate an incompressibility constraint, nonlinear advection, and two-way coupling between velocity and temperature. We consider the time-dependent Boussinesq equations at finite Rayleigh and Prandtl numbers. Our aim is to make high-order, fully implicit time integration practical by solving the resulting stage-coupled systems efficiently.

Projection and pressure-correction methods aim to reduce the solution cost by separating the velocity and pressure solves. They lead to relatively low cost per step, but their accuracy also depends on the splitting and its boundary treatment \cite{GuermondMinevShen2006}. Implicit--explicit (IMEX) methods \cite{ascher1997implicit} offer a different compromise by leaving selected terms explicit. Fully implicit Runge--Kutta methods instead retain the coupled equations at every stage. The Radau IIA family has classical order $2s-1$ with $s$ stages, is A-stable and L-stable, and is stiffly accurate \cite{HairerWanner1996}. These properties are attractive for stiff evolution equations with algebraic constraints. However, they do not eliminate accuracy restrictions on the time step or guarantee convergence of a nonlinear solver.

The efficient solution of the resulting linear systems remains computationally challenging. A dense Runge--Kutta matrix couples $s$ copies of the spatial problem, and the Newton systems inherit both this stage coupling and the saddle-point structure of incompressible flow. Treating the resulting matrix as a generic sparse system can make a high-order step expensive. This difficulty has motivated stage factorizations and block preconditioners \cite{PaznerPersson2017,RanaEtAl2021,SouthworthEtAl2021}, as well as monolithic multigrid methods that preserve stage coupling in the relaxation \cite{AbuLabdehEtAl2023}. Kirby~\cite{Kirby2024} relates convergence of a class of monolithic multigrid methods to associated single-stage problems with complex time steps. Such results motivate structured solvers, although they do not directly establish convergence for the nonlinear, thermally coupled systems studied here.

Preconditioning thermal convection also has a substantial history. Howle and Kirby~\cite{howle2012block} developed block preconditioners for incompressible flow with thermal convection. Subsequent work examined block triangular linearizations \cite{ke2017block}, augmented Lagrangian methods for steady buoyancy-driven flow \cite{KeAulisaDillonHowle2018}, multigrid for steady and unsteady convection \cite{KeAulisa2018}, and field-of-values estimates \cite{AulisaBorniaHowleKe2020}. The large-Prandtl asymptotic structure of the Boussinesq equations is studied in \cite{wang2007asymptotic}. In this work, we focus on time-dependent Rayleigh--B\'enard convection at finite Prandtl numbers using fully implicit Runge--Kutta methods.

We compare two approaches. The first adapts the augmented Lagrangian preconditioner approach of Leveque, He, and Olshanskii~\cite{LevequeHeOlshanskii2026} to the velocity--pressure subsystem and couples it multiplicatively to temperature. The second applies monolithic geometric multigrid to the full velocity--pressure--temperature stage system, with overlapping Vanka relaxation \cite{Vanka1986}. Both preconditioners aim to reduce the algebraic cost without compromising the accuracy of the fully implicit discretization. The augmented formulation also incorporates grad--div stabilization, resulting in a different discrete system from the standard formulation.

The spatial approximation provides another choice. Alongside Taylor--Hood elements \cite{TaylorHood1973}, we use the Alfeld--Sorokina pair in two dimensions \cite{AlfeldSorokina2016} and an enriched Guzm\'an--Neilan construction in three dimensions \cite{GuzmanNeilan2018}. Recent macroelement support in FIAT \cite{BrubeckKirby2025} makes these spaces available through Firedrake~\cite{RathgeberEtAl2016}. This allows us to assess the practical computational cost of improved discrete incompressibility in a coupled time-dependent problem.

The implementation uses Firedrake for spatial discretization and Irksome~\cite{FarrellKirbyMarchena2021} for fully implicit Runge–Kutta time integration. Irksome constructs the coupled Runge--Kutta forms from the semidiscrete variational problem; Firedrake's solver-composition interfaces \cite{KirbyMitchell2018} expose the field and stage structure to the preconditioners. Our contribution is a comprehensive assessment of these
fully implicit solvers for Rayleigh--B\'enard convection,
including spatial verification, mesh and parameter studies,
stage-dependent computational costs, comparison with a
pressure-correction method, and three-dimensional
solver performance and mass-conservation studies. Our results demonstrate that fully implicit Runge--Kutta methods can be computationally competitive for time-dependent Rayleigh--B\'enard convection when combined with effective multigrid preconditioners.

The remainder of the paper is organized as follows. \Cref{sec:model} introduces the incompressible Boussinesq equations and their variational formulation. \Cref{sec:discretization} describes the spatial finite element discretizations and fully implicit Runge--Kutta methods for temporal integration. \Cref{sec:algebraic} presents the resulting coupled nonlinear systems and their Newton linearization. \Cref{sec:augmented} discusses the augmented Lagrangian preconditioner formulation, and \Cref{sec:multigrid} describes the monolithic multigrid preconditioners with Vanka relaxation. \Cref{sec:numerical} presents the numerical experiments in two and three dimensions. Finally, we summarize our conclusions and future work.

\section{Rayleigh--B\'enard Convection}
\label{sec:model}

We consider the time-dependent Rayleigh--B\'enard convection
problem describing buoyancy-driven flow of an incompressible
fluid under the Boussinesq approximation. Let
$\Omega\subset\mathbb{R}^d$, $d=2,3$, be a bounded domain.
The nondimensional Boussinesq system of partial
differential equations (PDEs) is given by
\begin{subequations}\label{eq:boussinesq}
\begin{align}
\frac{\partial\mathbf{u}}{\partial t}
+(\mathbf{u}\cdot\nabla)\mathbf{u}
-\Pran\Delta\mathbf{u}
+\Pran\nabla p
&=-\Ra\Pran T\mathbf{g},
\label{eq:momentum}\\
\nabla\cdot\mathbf{u}&=0,
\label{eq:continuity}\\
\frac{\partial T}{\partial t}
+\mathbf{u}\cdot\nabla T-\Delta T&=0,
\label{eq:temperature}
\end{align}
\end{subequations}
where $\mathbf{u}$, $p$, and $T$ denote the velocity,
pressure, and temperature, respectively, and
$\mathbf{g}$ is the unit vector in the direction
of gravity. The Rayleigh and Prandtl numbers are
\begin{equation}\label{eq:parameters}
\Ra=\frac{g_0\alpha\Delta T H^3}{\nu\kappa},
\qquad
\Pran=\frac{\nu}{\kappa},
\end{equation}
where $\nu$ is the kinematic viscosity, $\kappa$ is the thermal diffusivity, $\alpha$ is the thermal expansion coefficient, $g_0$ is the magnitude of gravitational acceleration, and $H$ and $\Delta T$ denote the characteristic length and temperature difference, respectively.

The experiments are conducted on $\Omega=(0,1)^d$, with $y$ denoting the vertical coordinate and $\mathbf g=-\mathbf e_y$. Length, time, and velocity are scaled by $H$, $H^2/\kappa$, and $\kappa/H$, respectively. Pressure is scaled by $\rho_0\nu\kappa/H^2$, yielding the factor $\Pran$ in the pressure gradient term of \eqref{eq:momentum}.
We impose no-slip velocity boundary conditions,
fixed temperatures on the upper and lower boundaries,
and insulating lateral walls:
\begin{equation}\label{eq:boundary-conditions}
\begin{aligned}
\mathbf{u}&=\mathbf{0}
&&\text{on }\partial\Omega,\\
T&=1 &&\text{on }\Gamma_b,\\
T&=0 &&\text{on }\Gamma_t,\\
\nabla T\cdot\mathbf{n}&=0
&&\text{on }\Gamma_s,
\end{aligned}
\end{equation}
where $\Gamma_b$, $\Gamma_t$, and $\Gamma_s$
denote the bottom, top, and lateral boundaries,
respectively, and $\mathbf{n}$ is the outward
unit normal.

The initial velocity is zero, and the initial
temperature is a perturbation of the conductive
profile:
\begin{equation}\label{eq:initial-conditions}
\mathbf{u}_0=\mathbf{0},
\qquad
T_0(\mathbf x)=1-y+\varepsilon\cos(\pi x)\sin(\pi y),
\quad \varepsilon=0.05.
\end{equation}
The perturbation satisfies the prescribed thermal
boundary conditions and provides an initial
disturbance to the conductive state.

For the variational formulation, let
\begin{equation}\label{eq:function-spaces}
\begin{aligned}
V&=[H_0^1(\Omega)]^d,\qquad
Q=L_0^2(\Omega),\\
W&=\{S\in H^1(\Omega):
S=0\text{ on }\Gamma_b\cup\Gamma_t\},
\end{aligned}
\end{equation}
and let $W_D$ denote the affine temperature
space satisfying the prescribed Dirichlet
conditions. Here $L_0^2(\Omega)$ denotes the
space of square-integrable functions with
zero mean.

The weak problem is to find
$(\mathbf{u},p,T)\in V\times Q\times W_D$
such that, for all
$(\mathbf{v},q,S)\in V\times Q\times W$,
\begin{subequations}\label{eq:weak-boussinesq}
\begin{align}
\left(\frac{\partial\mathbf{u}}{\partial t},
\mathbf{v}\right)
+((\mathbf{u}\cdot\nabla)\mathbf{u},\mathbf{v})
+\Pran(\nabla\mathbf{u},\nabla\mathbf{v})
&\nonumber\\
-\Pran(p,\nabla\cdot\mathbf{v})
+\Ra\Pran(T\mathbf{g},\mathbf{v})
&=0,
\label{eq:weak-momentum}\\
(\nabla\cdot\mathbf{u},q)&=0,
\label{eq:weak-continuity}\\
\left(\frac{\partial T}{\partial t},S\right)
+(\mathbf{u}\cdot\nabla T,S)
+(\nabla T,\nabla S)&=0.
\label{eq:weak-temperature}
\end{align}
\end{subequations}
Here $(\cdot,\cdot)$ denotes the $L^2(\Omega)$
inner product, extended componentwise to vector-
and tensor-valued functions.

\section{Space and Time Discretization}
\label{sec:discretization}

\subsection{Spatial discretization}
\label{sec:spatial}

Let $\mathcal{T}_h$ be a conforming triangulation
of $\Omega$, with characteristic mesh size $h$.
We introduce finite-dimensional velocity,
pressure, and temperature spaces
\[
V_h\subset V,\qquad
Q_h\subset Q,\qquad
W_h\subset W,
\]
with $W_{D,h}$ denoting the corresponding affine
temperature trial space.

For the velocity--pressure pair, we consider continuous Taylor--Hood elements and Alfeld--Sorokina macroelements, with continuous piecewise linear elements for temperature. The Taylor--Hood pair consists of continuous piecewise quadratic velocities and continuous piecewise linear pressures. The Alfeld--Sorokina pair employs the corresponding macroelement spaces on Alfeld-refined meshes.

In two dimensions, the Alfeld--Sorokina velocity space consists of continuous piecewise quadratics on Alfeld-refined meshes with continuous divergence, paired with continuous piecewise linear pressures. The Alfeld–Sorokina pair is inf–sup stable in two dimensions, but the corresponding three-dimensional construction requires enrichment to achieve inf–sup stability. Following Guzm\'an and Neilan~\cite{GuzmanNeilan2018}, an inf--sup stable Stokes pair with continuous divergence is obtained by enriching the quadratic macroelement velocity space with cubic divergence-free face bubbles. In our three-dimensional experiments, we use the existing
FIAT implementation of this enriched element, available
in Firedrake under the name \texttt{GNH1div}. We denote this pair by GN in the three-dimensional tables and reserve AS for the two-dimensional pair. When $\nabla\cdot V_h\subseteq Q_h$, the discrete incompressibility constraint yields pointwise divergence-free velocities in exact arithmetic. In practice, small divergence errors may arise from algebraic solver tolerances. In contrast, Taylor--Hood elements enforce incompressibility only weakly and generally do not produce pointwise divergence-free velocities.

The semidiscrete Galerkin formulation is obtained
by restricting \eqref{eq:weak-boussinesq} to the
finite element spaces. It is convenient to write
this system in operator form:
\begin{subequations}\label{eq:semidiscrete}
\begin{align}
M_u\frac{d\mathbf{u}_h}{dt}
+\mathcal{N}_u(\mathbf{u}_h)
+\Pran K_u\mathbf{u}_h
-\Pran B^Tp_h
+\Ra\Pran G T_h
&=0,
\label{eq:semidiscrete-momentum}\\
B\mathbf{u}_h&=0,
\label{eq:semidiscrete-continuity}\\
M_T\frac{dT_h}{dt}
+\mathcal{N}_T(\mathbf{u}_h,T_h)
+K_T T_h
&=0.
\label{eq:semidiscrete-temperature}
\end{align}
\end{subequations}
Here $M_u$ and $M_T$ denote the velocity and temperature mass matrices, $K_u$ and $K_T$ the corresponding stiffness matrices, and $B$ the discrete divergence operator. The nonlinear operators $\mathcal{N}_u$ and $\mathcal{N}_T$ represent momentum and temperature advection, respectively. The operator $G$ denotes the temperature-to-momentum buoyancy coupling. In \eqref{eq:semidiscrete}, finite element functions are identified with their coefficient vectors, and contributions from prescribed boundary data are omitted for notational simplicity.
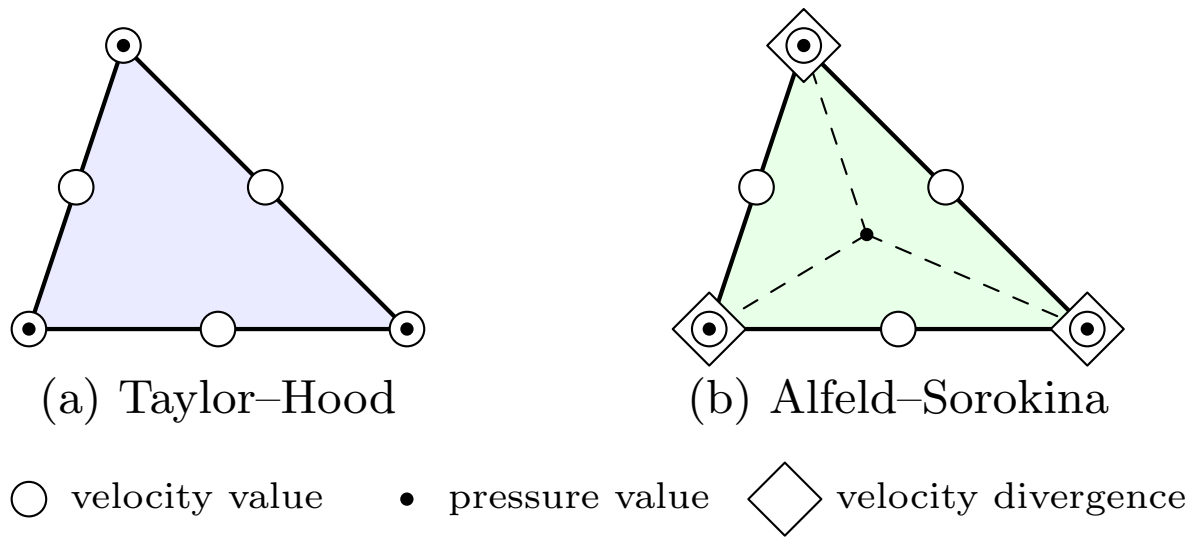
\begin{figure}[t]
\centering
\begin{tikzpicture}[
    scale=1.25,
    line cap=round,
    line join=round
]

\begin{scope}
\coordinate (A) at (0,0);
\coordinate (B) at (2,0);
\coordinate (C) at (0.5,1.5);

\fill[blue!8] (A)--(B)--(C)--cycle;
\draw[thick] (A)--(B)--(C)--cycle;

\foreach \P in {A,B,C}{
    \draw[fill=white] (\P) circle (2.6pt);
    \fill (\P) circle (1pt);
}

\foreach \X/\Y in {A/B,B/C,C/A}{
    \path (\X)--(\Y) coordinate[midway] (M);
    \draw[fill=white] (M) circle (2.6pt);
}

\node[font=\small] at (1,-0.38)
    {(a) Taylor--Hood};
\end{scope}

\begin{scope}[xshift=3.6cm]
\coordinate (A) at (0,0);
\coordinate (B) at (2,0);
\coordinate (C) at (0.5,1.5);
\coordinate (G) at ({2.5/3},0.5);

\fill[green!9] (A)--(B)--(C)--cycle;
\draw[thick] (A)--(B)--(C)--cycle;

\draw[dashed] (G)--(A);
\draw[dashed] (G)--(B);
\draw[dashed] (G)--(C);

\foreach \P in {A,B,C}{
    \begin{scope}[shift={(\P)}]
        \draw[fill=white]
            (0,5.5pt)--(5.5pt,0)--
            (0,-5.5pt)--(-5.5pt,0)--cycle;
        \draw[fill=white] (0,0) circle (2.6pt);
        \fill (0,0) circle (1pt);
    \end{scope}
}

\foreach \X/\Y in {A/B,B/C,C/A}{
    \path (\X)--(\Y) coordinate[midway] (M);
    \draw[fill=white] (M) circle (2.6pt);
}

\fill (G) circle (1pt);

\node[font=\small] at (1,-0.38)
    {(b) Alfeld--Sorokina};
\end{scope}

\begin{scope}[yshift=-0.9cm]
    \draw[fill=white] (0,0) circle (2.6pt);
    \node[right,font=\scriptsize] at (0.12,0)
        {velocity value};

    \fill (2,0) circle (1pt);
    \node[right,font=\scriptsize] at (2.12,0)
        {pressure value};

    \begin{scope}[shift={(4,0)}]
        \draw[fill=white]
            (0,5.5pt)--(5.5pt,0)--
            (0,-5.5pt)--(-5.5pt,0)--cycle;
    \end{scope}
    \node[right,font=\scriptsize] at (4.17,0)
        {velocity divergence};
\end{scope}

\end{tikzpicture}

\caption{Degrees of freedom for the two-dimensional
velocity--pressure pairs.
(a) Taylor--Hood $P_2/P_1$.
(b) The quadratic Alfeld--Sorokina velocity macroelement
with continuous divergence, paired with continuous
piecewise-linear pressure on the Alfeld split.
Each open circle represents two velocity components,
each dot a pressure value, and each diamond a
velocity-divergence value. Dashed lines mark the internal
edges of the split.}
\label{fig:fe-comparison}
\end{figure}
\subsection{Fully implicit Runge--Kutta discretization}
\label{sec:time}

For temporal discretization, we employ an
$s$-stage fully implicit Runge--Kutta method
with Butcher tableau
\begin{equation}\label{eq:butcher}
\begin{array}{c|c}
\mathbf{c} & A_{\mathrm{RK}}\\
\hline
\noalign{\vskip 3pt}
&\mathbf{b}^T
\end{array},
\qquad
A_{\mathrm{RK}}=(a_{ij})_{i,j=1}^s.
\end{equation}
Here $\mathbf{c}$ contains the stage nodes,
$\mathbf{b}$ the quadrature weights, and
$A_{\mathrm{RK}}$ the Runge--Kutta coefficient
matrix. We focus on the Radau IIA family,
which is A-stable and L-stable, and has
classical order $2s-1$ for $s$ stages~\cite{HairerWanner1996}.

Let $t_n$ denote the time at step $n$, where $n$ is the time-step index, and let $\Delta t=t_{n+1}-t_n$ denote the time-step size. For each stage
$i=1,\ldots,s$, let
$(\mathbf{U}_i,P_i,\Theta_i)$ denote the
velocity, pressure, and temperature stage
approximations. Introducing the stage
derivatives $\mathbf{k}_i$ and $\ell_i$,
the Runge--Kutta relations are
\begin{equation}\label{eq:stage-relations}
\begin{aligned}
\mathbf{U}_i
&=\mathbf{u}_h^n
+\Delta t\sum_{j=1}^s a_{ij}\mathbf{k}_j,\\
\Theta_i
&=T_h^n
+\Delta t\sum_{j=1}^s a_{ij}\ell_j.
\end{aligned}
\end{equation}
At each stage, the semidiscrete equations
\eqref{eq:semidiscrete} become
\begin{subequations}\label{eq:irk-stages}
\begin{align}
M_u\mathbf{k}_i
+\mathcal{N}_u(\mathbf{U}_i)
+\Pran K_u\mathbf{U}_i
-\Pran B^TP_i
+\Ra\Pran G\Theta_i
&=0,
\label{eq:irk-momentum}\\
B\mathbf{U}_i&=0,
\label{eq:irk-continuity}\\
M_T\ell_i
+\mathcal{N}_T(\mathbf{U}_i,\Theta_i)
+K_T\Theta_i
&=0.
\label{eq:irk-temperature}
\end{align}
\end{subequations}
The Runge--Kutta stages are coupled through the generally dense coefficient matrix $A_{\mathrm{RK}}$ and are solved simultaneously.

After computing the stages, the velocity and
temperature are updated according to
\begin{equation}\label{eq:irk-update}
\begin{aligned}
\mathbf{u}_h^{n+1}
&=\mathbf{u}_h^n
+\Delta t\sum_{i=1}^s b_i\mathbf{k}_i,\\
T_h^{n+1}
&=T_h^n
+\Delta t\sum_{i=1}^s b_i\ell_i.
\end{aligned}
\end{equation}
The pressure is an algebraic variable determined
by the incompressibility constraint at each stage.

Radau IIA satisfies $c_s=1$ and $b_j=a_{sj}$, so the last stage gives the new velocity and temperature, and we take $p_h^{n+1}=P_s$. Pressure constants are fixed separately at every stage. For example, the two-stage method is
\[
\begin{array}{c|cc}
1/3&5/12&-1/12\\
1&3/4&1/4\\\hline
&3/4&1/4
\end{array}.
\]
The three- and four-stage Radau IIA methods have classical orders five and seven, respectively. However, spatial discretization errors, boundary conditions, and algebraic solver tolerances may affect the observed temporal accuracy. We therefore evaluate temporal accuracy using a refined numerical reference solution.

\section{Stage-coupled algebraic systems}\label{sec:algebraic}
Since the Runge--Kutta coefficient matrix $A_{\mathrm{RK}}$ is invertible, the stage derivatives in \eqref{eq:stage-relations} can be eliminated:
\begin{equation}\label{eq:stage-derivative-inverse}
\begin{aligned}
\mathbf{k}_i
&=\frac{1}{\Delta t}
\sum_{j=1}^s (A_{\mathrm{RK}}^{-1})_{ij}
(\mathbf{U}_j-\mathbf{u}_h^n),\\
\ell_i
&=\frac{1}{\Delta t}
\sum_{j=1}^s (A_{\mathrm{RK}}^{-1})_{ij}
(\Theta_j-T_h^n).
\end{aligned}
\end{equation}
Let $\mathbf U=(\mathbf U_1,\ldots,\mathbf U_s)^T$, $\mathbf P=(P_1,\ldots,P_s)^T$, and $\boldsymbol\Theta=(\Theta_1,\ldots,\Theta_s)^T$ collect the stage variables. Writing $x=(\mathbf U,\mathbf P,\boldsymbol\Theta)^T$, Newton's method applied to the coupled stage residual $\mathcal R(x)$ yields
\begin{equation}\label{eq:newton-system}
\mathcal J_m\delta x_m=-\mathcal R(x_m),
\end{equation}
where $\mathcal J_m=\mathcal R'(x_m)$ is the Jacobian evaluated at the current Newton iterate $x_m$. The linear system is solved inexactly, with the correction $\delta x_m$ satisfying
\begin{equation}\label{eq:inexact-newton}
\|\mathcal J_m\delta x_m+\mathcal R(x_m)\|
\leq\eta_m\|\mathcal R(x_m)\|,
\end{equation}
subject to the prescribed absolute linear tolerance. In our solver implementation, we employ Eisenstat--Walker forcing terms~\cite{EisenstatWalker1996} to adapt the linear solver tolerances during Newton iteration. We use an $L^2$ residual-norm line search to improve nonlinear convergence and flexible GMRES~\cite{Saad1993} to accommodate iterative preconditioner applications. Solver performance is assessed through iteration counts and computational cost.

Ordering the unknowns by physical field, the Jacobian takes the form
\begin{equation}\label{eq:coupled-jacobian}
\mathcal{J}=
\begin{bmatrix}
\mathcal{F} & -\Pran\mathcal{B}^T
            & \mathcal{C}_{uT}\\
\mathcal{B} & 0 & 0\\
\mathcal{C}_{Tu} & 0 & \mathcal{H}
\end{bmatrix},
\end{equation}
where $\mathcal{B}=I_s\otimes B$ denotes the stage-wise divergence operator. The blocks $\mathcal{F}$ and $\mathcal{H}$ correspond to the velocity and temperature equations, respectively, and are given by
\begin{equation}\label{eq:stage-diagonal-blocks}
\begin{aligned}
\mathcal{F}
&=\frac{1}{\Delta t}
A_{\mathrm{RK}}^{-1}\otimes M_u
+\operatorname{diag}(F_1,\ldots,F_s),\\
\mathcal{H}
&=\frac{1}{\Delta t}
A_{\mathrm{RK}}^{-1}\otimes M_T
+\operatorname{diag}(H_1,\ldots,H_s),
\end{aligned}
\end{equation}
with
\[
F_i=\mathcal{N}_u'(\mathbf{U}_i)
+\Pran K_u,
\qquad
H_i=\mathcal{N}_{T,T}'(\mathbf{U}_i,\Theta_i)
+K_T.
\]
Here the primes denote the corresponding Fr\'echet derivatives. The off-diagonal blocks $\mathcal{C}_{uT}$ and $\mathcal{C}_{Tu}$ represent the linearized buoyancy force and the dependence of thermal advection on velocity, respectively. Both are block diagonal with respect to the stage index.

The Jacobian \eqref{eq:coupled-jacobian} combines the saddle-point structure of incompressible flow with velocity--temperature coupling. The terms involving $A_{\mathrm{RK}}^{-1}$ introduce additional coupling between the Runge--Kutta stages, requiring the solution of a coupled multistage system.

We consider two preconditioning strategies for \eqref{eq:newton-system}. The first employs an augmented Lagrangian formulation with a block decomposition separating the velocity--pressure subsystem from temperature. The second applies geometric multigrid directly to the coupled stage system.

\section{Augmented Lagrangian Preconditioning}
\label{sec:augmented}

Our first preconditioning approach extends the augmented Lagrangian
preconditioner for fully implicit Runge--Kutta discretizations of the
incompressible Navier--Stokes equations~\cite{LevequeHeOlshanskii2026}
to the coupled Boussinesq system. The approach incorporates grad--div
stabilization into the momentum equation and employs a block
preconditioner that exploits the velocity--pressure structure and
accounts for the coupling with temperature.

We augment the momentum variational form with
\begin{equation}\label{eq:grad-div}
\gamma(\nabla\cdot\mathbf{u}_h,
       \nabla\cdot\mathbf{v}_h),
\qquad \gamma>0,
\end{equation}
where $\gamma$ is the augmentation parameter. For discretizations
that enforce incompressibility pointwise, the grad--div term
vanishes at the discrete solution, so the augmentation preserves
that solution and modifies the velocity block used in the solve.
For Taylor--Hood elements, incompressibility is enforced weakly,
and the grad--div term generally changes the discrete solution
by penalizing the velocity divergence. In our implementation,
the augmentation is incorporated directly into the variational
formulation before applying the fully implicit Runge--Kutta
discretization.

Let $D$ denote the matrix associated with the grad--div bilinear
form. In the stage-value representation, the augmented velocity
block becomes
\begin{equation}\label{eq:augmented-velocity}
\mathcal{F}_{\gamma}
=
\mathcal{F}+\gamma(I_s\otimes D).
\end{equation}
The corresponding velocity--pressure subsystem is
\begin{equation}\label{eq:augmented-up}
\mathcal{J}_{up,\gamma}
=
\begin{bmatrix}
\mathcal{F}_{\gamma} & -\Pran\mathcal{B}^T\\
\mathcal{B} & 0
\end{bmatrix}.
\end{equation}
The augmented formulation therefore preserves the Runge--Kutta
stage coupling and modifies the velocity operator. We note
that $D$ is obtained from the variational grad--div term and
is not generally identical to the projected algebraic operator
$B^TM_p^{-1}B$.

\subsection{Block preconditioning}
\label{sec:al-block}

We first describe the block decomposition used to precondition
the complete Boussinesq system. Grouping the velocity and
pressure unknowns together, the augmented Jacobian can be
written as
\begin{equation}\label{eq:outer-al-block}
\mathcal{J}_{\gamma}
=
\begin{bmatrix}
\mathcal{J}_{up,\gamma}
& \mathcal{C}_{up,T}\\
\mathcal{C}_{T,up} & \mathcal{H}
\end{bmatrix},
\end{equation}
where $\mathcal{C}_{up,T}$ and $\mathcal{C}_{T,up}$ represent
the coupling between the velocity--pressure and temperature
subsystems.

We employ a multiplicative field-split preconditioner based on
this partition. The velocity--pressure subsystem is treated
using an upper triangular Schur complement preconditioner. The temperature block is approximately inverted using
algebraic multigrid. All algebraic multigrid solvers in the augmented
Lagrangian preconditioner are implemented using hypre's BoomerAMG~\cite{falgout2002hypre,HensonYang2002}. This construction exploits the
velocity--pressure saddle-point structure and incorporates
the coupling with temperature through the outer block
decomposition.

The augmented velocity block $\mathcal{F}_\gamma$ is
approximately inverted using a multiplicative splitting
across the Runge--Kutta stages. Each stage block is treated
using ten GMRES iterations with an assembled algebraic
multigrid preconditioner. This provides a stage-block
Gauss--Seidel approximation within the coupled
Runge--Kutta system.

The remaining component of the velocity--pressure
preconditioner is an approximation to the pressure Schur
complement, which we describe as follows.

\subsection{Schur complement approximation}
\label{sec:schur}

For the velocity--pressure system \eqref{eq:augmented-up},
the pressure Schur complement is
\begin{equation}
\mathcal{S}_{\gamma}
=
\Pran\mathcal{B}
\mathcal{F}_{\gamma}^{-1}
\mathcal{B}^T.
\end{equation}
Direct application of its inverse would require solving
systems involving the augmented velocity operator.
We therefore construct an approximation based on the
pressure mass and stiffness matrices, following the
approach of~\cite{LevequeHeOlshanskii2026}.

Let $M_p$ and $K_p$ denote the pressure mass and
stiffness matrices, respectively. In the stage-value
representation, a mass-plus-stiffness approximation
to the inverse Schur complement is
\begin{equation}
\widetilde{\mathcal{S}}_{\gamma}^{-1}
=
\frac{1}{\Pran\Delta t}
\left(
A_{\mathrm{RK}}^{-1}\otimes K_p^{-1}
\right)
+
\frac{\Pran+\gamma}{\Pran}
\left(
I_s\otimes M_p^{-1}
\right).
\end{equation}
The first term is associated with the time-dependent
contribution, and the second reflects the viscosity
and grad--div augmentation. This approximation is
motivated by the unsteady Stokes operator; convection
and boundary effects are not represented exactly.

To obtain the corresponding approximation in the
stage-derivative formulation, we introduce the
pressure normalization $\widehat p=\Pran p$ and define
\[
\widehat{\mathcal{S}}_{\gamma,d}
=
\frac{\mathcal{S}_{\gamma}}{\Pran}
(\Delta t A_{\mathrm{RK}}\otimes I_p).
\]
The transformation modifies the stage factors appearing
in the inverse Schur approximation.

To obtain the form used in the preconditioner, define
\[
\mathcal{W}
=
\Delta t A_{\mathrm{RK}}\otimes M_p.
\]
We consider the additive approximation
\begin{equation}\label{eq:schur-initial}
\widetilde{\widehat{\mathcal{S}}}_{\gamma,d}^{-1}
=
\gamma\mathcal{W}^{-1}
+
\frac{1}{\Delta t^2}
(A_{\mathrm{RK}}^{-1}\otimes I_p)
\widetilde S_{\mathrm{int}}^{-1}
(A_{\mathrm{RK}}^{-1}\otimes I_p),
\end{equation}
where
\begin{equation}\label{eq:schur-intermediate}
\widetilde S_{\mathrm{int}}^{-1}
=
I_s\otimes K_p^{-1}
+
\Pran\Delta t\,
A_{\mathrm{RK}}\otimes M_p^{-1}.
\end{equation}
Using the mixed-product property of Kronecker products
and the identity
\[
\mathcal{W}^{-1}
=
\frac{1}{\Delta t}
A_{\mathrm{RK}}^{-1}\otimes M_p^{-1},
\]
we obtain the two-term approximation
\begin{equation}\label{eq:schur-kron}
\widetilde{\widehat{\mathcal{S}}}_{\gamma,d}^{-1}
={}
\frac{\Pran+\gamma}{\Delta t}
\left(
A_{\mathrm{RK}}^{-1}\otimes M_p^{-1}
\right)\\
+
\frac{1}{\Delta t^2}
\left(
A_{\mathrm{RK}}^{-2}\otimes K_p^{-1}
\right).
\end{equation}

The approximation \eqref{eq:schur-kron} consists of
pressure mass and stiffness contributions, each coupled
across stages through a small dense Runge--Kutta matrix.
The pressure inverses are understood on the appropriate
pressure space, with the constant-pressure nullspace
handled separately.

In our implementation, the two terms in
\eqref{eq:schur-kron} are combined additively using
custom \texttt{MassKronPC} and \texttt{StiffnessKronPC} \footnote{The implementation is available in
the \texttt{436ahsan/kron-test} branch of Irksome:
\url{https://github.com/firedrakeproject/Irksome/tree/436ahsan/kron-test}.}
components within Firedrake and Irksome. The corresponding
stage matrices are
\[
L_M
=
\frac{\Pran+\gamma}{\Delta t}
A_{\mathrm{RK}}^{-1},
\qquad
L_K
=
\frac{1}{\Delta t^2}
A_{\mathrm{RK}}^{-2},
\]
respectively. Algebraic multigrid is used to
approximate the spatial pressure inverses.

For the Taylor--Hood experiments, we use the augmentation
parameter $\gamma=1$. The complete nested preconditioner
is applied within flexible GMRES to solve the Newton
linear systems.

\section{Monolithic Multigrid}
\label{sec:multigrid}

Our second preconditioning strategy employs monolithic geometric
multigrid for the coupled velocity--pressure--temperature system.
Unlike the augmented Lagrangian approach, which relies on block factorization and pressure Schur complement approximations, the monolithic method applies multigrid directly to the complete stage-coupled Jacobian. This approach builds on monolithic multigrid methods
for fully implicit Runge--Kutta discretizations of incompressible flow~\cite{Kirby2024}.

Let
\[
\mathcal{T}_0\prec\mathcal{T}_1
\prec\cdots\prec\mathcal{T}_L
\]
denote a sequence of successively refined meshes, with
corresponding mixed finite element spaces
$\mathcal{X}_\ell$.
For an $s$-stage Runge--Kutta method, the discrete system contains $s$ copies of the velocity, pressure, and temperature unknowns. We construct a geometric multigrid hierarchy for the complete coupled system, using finite element transfer operators for prolongation and restriction between levels.

An essential component of the multigrid preconditioner
is an effective relaxation method for the coupled
saddle-point system. We employ overlapping additive
Schwarz relaxation based on vertex-centered Vanka patches, implemented using \texttt{ASMVankaPC} in
Firedrake~\cite{RathgeberEtAl2016}. These patches provide
local coupling between velocity and pressure, together
with the buoyancy and thermal-advection contributions
arising from the Boussinesq equations.

For each mesh vertex $v$, we construct a patch containing
the velocity degrees of freedom associated with the
closure of its cell star across all Runge--Kutta stages.
The pressure and temperature degrees of freedom are
restricted to the entities in the vertex star.
The \texttt{pc\_vanka\_exclude\_subspaces} option
controls this restriction without removing pressure
or temperature from the local patch systems.
Consequently, the patch solves retain the coupling
between the physical variables and Runge--Kutta stages.

Let $J_\ell$ denote the Jacobian on multigrid level
$\ell$, and let $R_v$ be the restriction operator
associated with a vertex patch. The local patch matrix
is $J_v=R_vJ_\ell R_v^T$, and the additive Schwarz
preconditioner is given by
\begin{equation}\label{eq:patch-action}
\mathcal M_\ell^{-1}r
=
\sum_v R_v^T J_v^{-1}R_v r.
\end{equation}
The local systems are solved using sparse LU
factorization. By restricting the
pressure and temperature degrees of freedom,
the patches remain smaller than those obtained
from the full closure without sacrificing the local
saddle-point and thermal coupling. This extends
Vanka-type relaxation~\cite{Vanka1986} to the
coupled multistage Boussinesq system.

On each multigrid level except the coarsest, we apply ten
GMRES iterations preconditioned by the additive
Schwarz operator \eqref{eq:patch-action}. The
coarse-grid system is solved directly using
MUMPS \cite{AmestoyEtAl2001}, and the resulting multigrid preconditioner
is applied within outer FGMRES.

We employ the same multigrid framework for
Taylor--Hood and macroelement discretizations.
The effectiveness and computational cost of
the monolithic and augmented Lagrangian
preconditioners are examined in
the next section.

\section{Numerical results}\label{sec:numerical}\label{sec:numerics}
In this section, we investigate the accuracy and computational performance of the proposed solvers for two- and three-dimensional Rayleigh--B\'enard convection. We first verify spatial convergence using a manufactured solution, followed by studies of mesh refinement, Rayleigh and Prandtl numbers, and the number of Runge--Kutta stages. We then examine parallel performance, solver reuse, and the tradeoff between temporal accuracy and computational cost, including a comparison with incremental pressure
correction. Finally, we present three-dimensional experiments to assess solver robustness and parallel scalability.

\subsection{Experimental setting}
The two-dimensional benchmark is posed on $\Omega=(0,1)^2$, with no-slip velocity conditions on all boundaries. The temperature is fixed at $T=1$ on the bottom boundary and $T=0$ on the top boundary, with insulating conditions on the vertical walls. The fluid is initially at rest,
and the initial temperature is prescribed by
\[
T_0(x,y)=1-y+0.05\cos(\pi x)\sin(\pi y).
\]
The simulations are carried out up to the final time $t_f=2$. Uniform refinement of a $16\times16$ base mesh yields resolutions $N=64,128,256$, with additional experiments at $N=512$. Time integration uses two-, three-, and four-stage Radau IIA methods. Depending on the experiment, we choose $\Delta t=1/N$ or hold the time step fixed at $\Delta t=1/128$, as specified in each study.

We denote augmented Lagrangian preconditioning with Taylor--Hood elements by AL--TH, and monolithic multigrid with Taylor--Hood and Alfeld--Sorokina elements by MG--TH and MG--AS, respectively. Both discretizations use continuous piecewise linear elements for temperature. At $N=128$, the mixed finite element spaces have $165{,}380$ and $214{,}789$ degrees of freedom for Taylor--Hood and Alfeld--Sorokina, respectively. At $N=256$, these dimensions increase to $658{,}436$ and $855{,}557$. The corresponding $s$-stage systems contain $s$ times as many unknowns. Thus, comparisons at the same mesh resolution involve different algebraic problem sizes. Grad–div stabilization with $\gamma=1$ is included only in the AL–TH configuration.

The numerical implementation uses Firedrake and Irksome~\cite{RathgeberEtAl2016,FarrellKirbyMarchena2021}. For the monolithic experiments, Newton's method employs an $L^2$ residual-norm line search with relative and absolute tolerances of $10^{-5}$. The linearized systems are solved using FGMRES with a maximum of 200 iterations, and an absolute tolerance of $10^{-5}$. The relative tolerance is adjusted using Eisenstat--Walker forcing terms~\cite{EisenstatWalker1996}. On each noncoarse multigrid level, we apply ten GMRES iterations preconditioned by vertex-centered Vanka relaxation. The local patch systems are solved using UMFPACK~\cite{Davis2004}, and the coarse-grid system is solved using MUMPS~\cite{AmestoyEtAl2001}. Unless otherwise specified, the Jacobian and preconditioner are recomputed at each Newton iteration.

Let $n_t$, $n_{\mathrm{N}}$, and $n_{\mathrm{L}}$ denote the
total numbers of time steps, Newton iterations, and outer
Krylov iterations, respectively. We report the average
Newton iterations per time step, outer Krylov iterations
per Newton iteration, and computational time per step:
\begin{equation}\label{eq:performance-metrics}
\overline n_{\mathrm{N}}=\frac{n_{\mathrm{N}}}{n_t},
\qquad
\overline n_{\mathrm{L}}=\frac{n_{\mathrm{L}}}{n_{\mathrm{N}}},
\qquad
\overline t=\frac{t_{\mathrm{solve}}}{n_t}.
\end{equation}
Values of $\overline n_{\mathrm{N}}<1$ arise when some time
steps satisfy the nonlinear convergence criterion without
requiring a Newton update. The reported Krylov iteration
counts exclude inner multigrid and block-solver iterations; therefore, we consider both iteration counts and elapsed time when evaluating solver performance.

All two-dimensional experiments were performed on a Mac Studio equipped with an Apple M5 Max processor (18-core CPU, 40-core GPU) and 128 GB of memory. Unless otherwise specified, we used 16 MPI processes with one OpenMP thread per process. All computations were performed on the CPU. The three-dimensional experiments were conducted on the Magnolia cluster at the Mississippi Center for Supercomputing Research (MCSR), with MPI process counts specified for each experiment. Reported timings are elapsed wall-clock times from individual runs, with visualization disabled during time stepping. 

\begin{figure}[!htbp]\centering
\includegraphics[width=.96\linewidth]{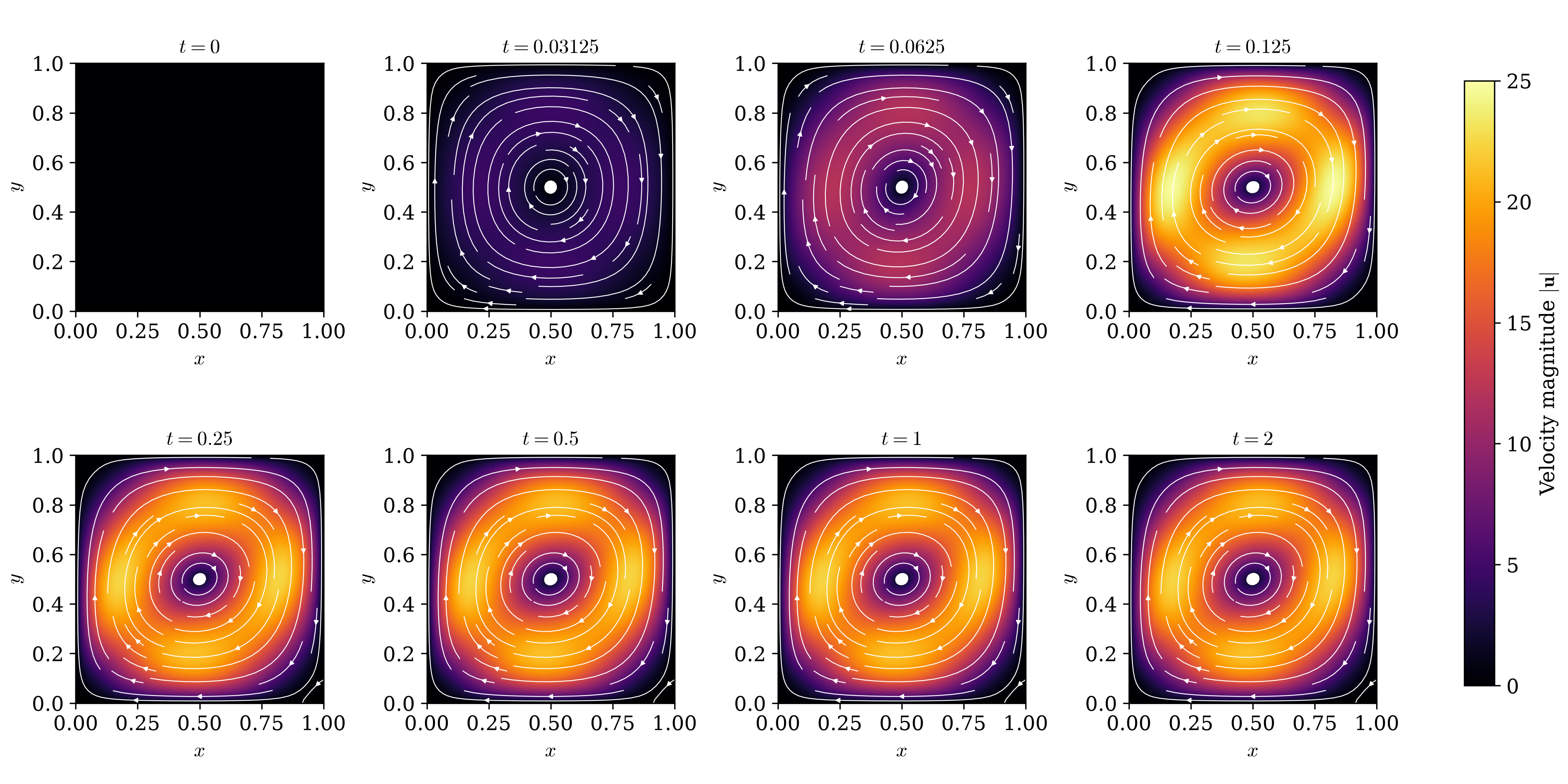}\\[-2mm]
\includegraphics[width=.96\linewidth]{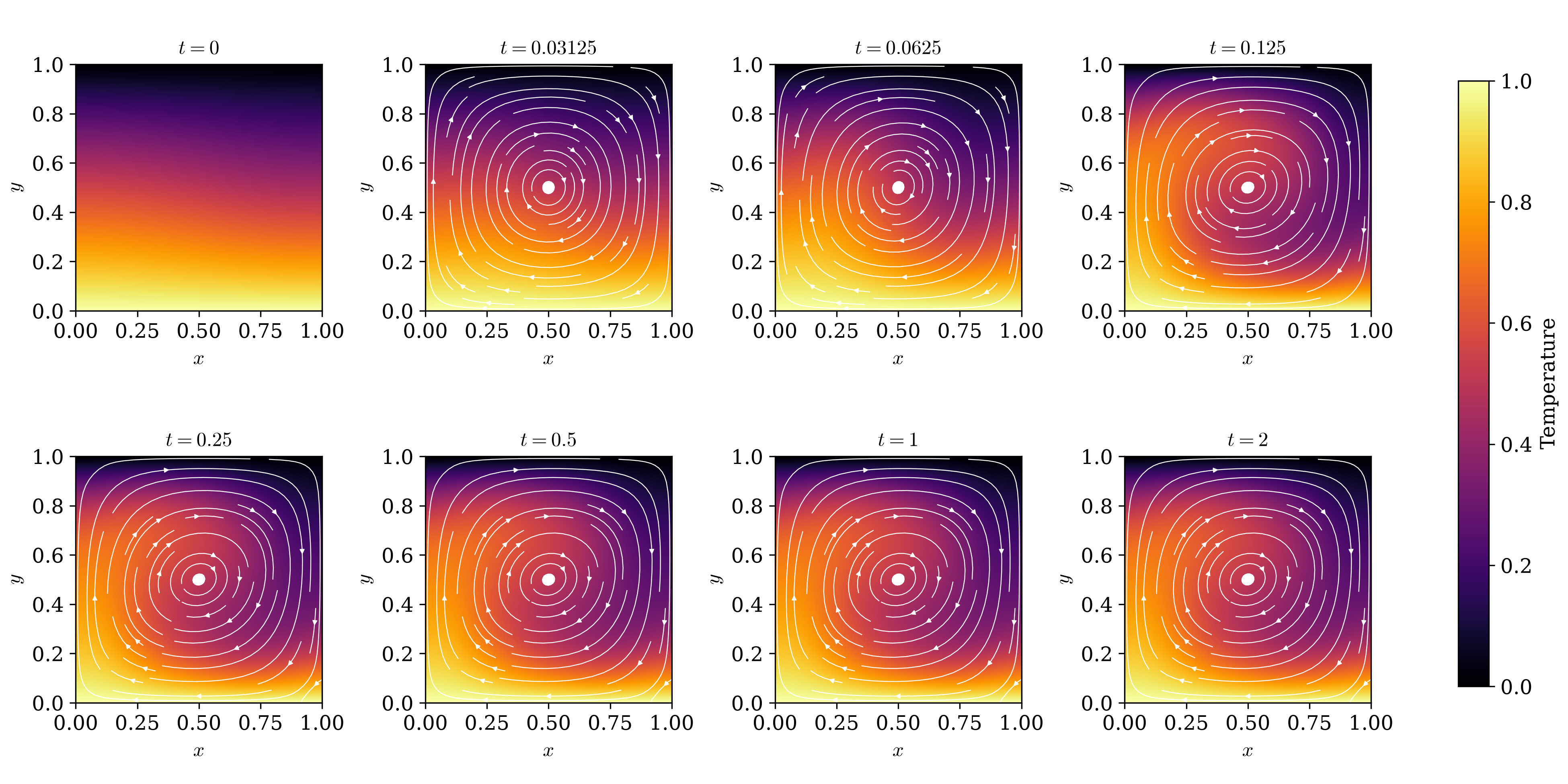}
\caption{Two-dimensional benchmark snapshots at the indicated times, showing velocity magnitude with streamlines (top two rows) and temperature with streamlines (bottom two rows). The snapshots illustrate the development and persistence of the convection roll.}\label{fig:two-d-snapshots}
\end{figure}

Figure~\ref{fig:two-d-snapshots} illustrates the evolution of the velocity and temperature fields. Starting from rest, the flow develops a circulation that strengthens during the initial transient and forms a persistent convection roll. Meanwhile, thermal advection progressively deforms the temperature contours. The most significant changes occur during the early transient, motivating comparisons of numerical solutions at multiple times rather than only at the final time.

\subsection{Spatial accuracy}
We first verify the spatial accuracy of the finite element
discretizations using the method of manufactured solutions.
For $\Ra=\Pran=1$, we prescribe the exact solution
\begin{equation}\label{eq:numerical-manufactured}
\begin{aligned}
\psi(x,y)&=\sin^2(\pi x)\sin^2(\pi y),
\qquad a(t)=e^{4t},\\
\mathbf u_e&=a(t)(\partial_y\psi,-\partial_x\psi)^T,\\
p_e&=a(t)\cos(\pi x)\sin(\pi y),\\
T_e&=1-y+0.05a(t)\cos(\pi x)\sin(\pi y).
\end{aligned}
\end{equation}
The corresponding forcing terms are obtained by substituting
these expressions into the governing equations. The prescribed
velocity is divergence-free and satisfies the no-slip boundary
conditions. The temperature satisfies the prescribed thermal boundary conditions.

We employ the three-stage Radau IIA method with
$\Delta t=5\times10^{-4}$ and final time $t_f=0.01$.
The Newton linear systems are solved using a sparse direct
solver, with relative and absolute nonlinear tolerances of
$10^{-9}$ for Taylor--Hood and $10^{-7}$ for Alfeld--Sorokina.
For the Alfeld--Sorokina discretization, the initial velocity
is obtained by $L^2$ projection of the exact initial condition.

\begin{table}[!htbp]
\centering\small
\caption{Spatial convergence for the manufactured solution at $t_f=0.01$. The quantities $E_u$, $E_p$, and $E_T$ denote the $L^2$ errors in velocity, pressure, and temperature, respectively, with the pressure error computed modulo constants. The columns $r_u$ and $r_T$ report the observed convergence orders under successive mesh refinement.}\label{tab:spatial-results}
\setlength{\tabcolsep}{5pt}
\begin{tabular}{lrrrrrr}
\toprule
Pair & $N$ & $E_u$ & $r_u$ & $E_p$ & $E_T$ & $r_T$\\\midrule
TH & 16 & $1.363\times10^{-3}$ & --- & $3.021\times10^{-3}$ & $2.143\times10^{-4}$ & ---\\
TH & 32 & $1.733\times10^{-4}$ & 2.98 & $4.646\times10^{-4}$ & $5.373\times10^{-5}$ & 2.00\\
TH & 64 & $2.176\times10^{-5}$ & 2.99 & $1.059\times10^{-4}$ & $1.344\times10^{-5}$ & 2.00\\
AS & 16 & $3.007\times10^{-3}$ & --- & $7.195\times10^{-2}$ & $2.123\times10^{-4}$ & ---\\
AS & 32 & $3.811\times10^{-4}$ & 2.98 & $1.084\times10^{-2}$ & $5.358\times10^{-5}$ & 1.99\\
AS & 64 & $4.756\times10^{-5}$ & 3.00 & $1.579\times10^{-3}$ & $1.343\times10^{-5}$ & 2.00\\
\bottomrule\end{tabular}
\end{table}

\begin{figure}[!htbp]\centering
\includegraphics[width=\linewidth]{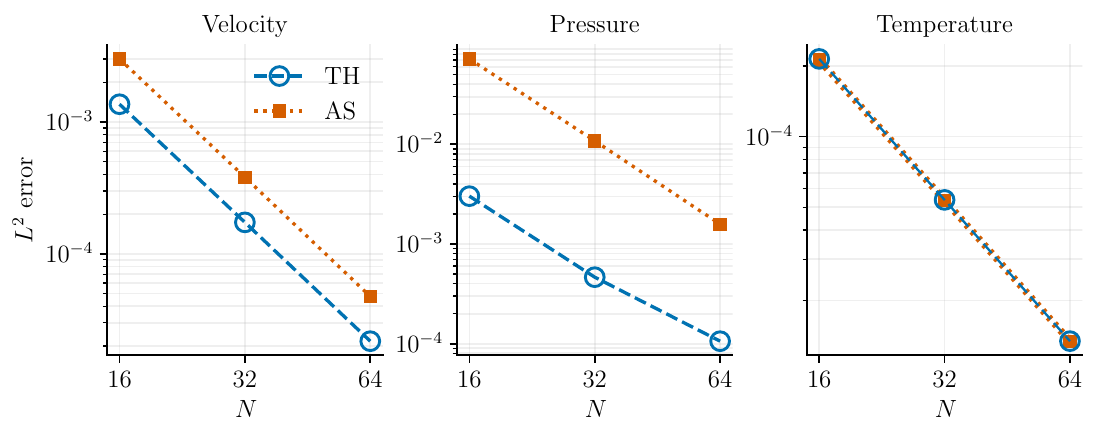}
\caption{Spatial errors for the manufactured solution. Both element pairs approach third-order velocity and second-order temperature convergence.}\label{fig:spatial-results}
\end{figure}

Table~\ref{tab:spatial-results} and Figure~\ref{fig:spatial-results} demonstrate approximately third-order convergence in velocity and second-order convergence in temperature for both finite element pairs. The pressure errors also decrease under mesh refinement, with observed convergence rates of $2.13$ for Taylor--Hood and $2.78$ for Alfeld--Sorokina on the finest refinement. To assess the influence of temporal discretization, we halve the time step at $N=64$ and observe no change in the reported errors to the tabulated precision. These results confirm the expected spatial convergence behavior for velocity and temperature in the manufactured solution test.

\subsection{Mesh and physical-parameter dependence}
\begin{figure}[!htbp]\centering
\includegraphics[width=\linewidth]{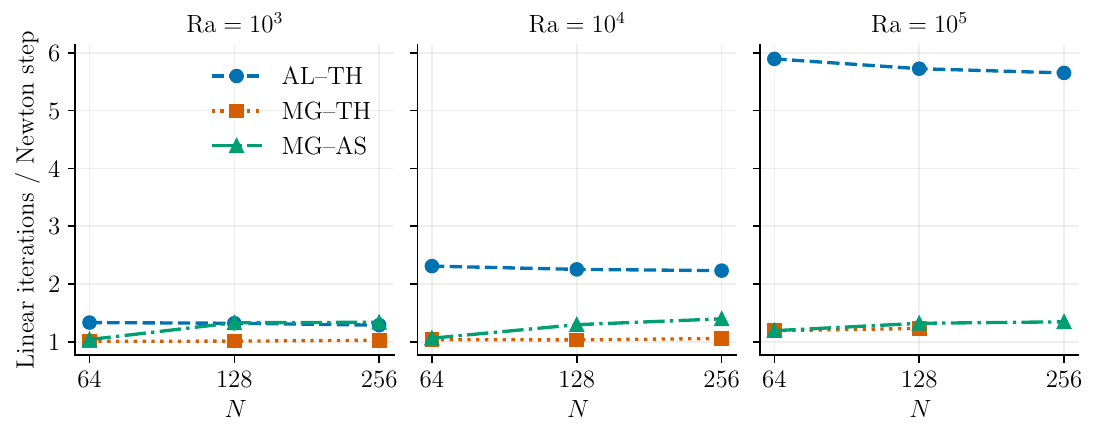}
\caption{Average outer Krylov iterations per Newton iteration
versus mesh resolution for $\Ra=10^3,10^4,10^5$, with
$\Delta t=1/128$, $s=2$, $\Pran=1$, and 16 MPI processes.
All reported results correspond to completed integrations.
The MG--TH simulation at $\Ra=10^5$, $N=256$ did not converge and is therefore omitted.}\label{fig:mesh-results}
\end{figure}

We first investigate the effect of mesh refinement by fixing $\Delta t=1/128$, $s=2$, and $\Pran=1$. Figure~\ref{fig:mesh-results} presents the average outer Krylov iterations per Newton iteration for $\Ra=10^3,10^4,10^5$. The AL--TH preconditioner exhibits nearly mesh-independent convergence across all three Rayleigh numbers. For example, at $\Ra=10^4$, the average iteration count decreases slightly from $2.31$ at $N=64$ to $2.23$ at $N=256$. However, the iteration counts increase with $\Ra$, from approximately $1.3$ at $\Ra=10^3$ to $5.7$--$5.9$ at $\Ra=10^5$.

The monolithic multigrid method exhibits consistently low outer Krylov iteration counts with both finite element discretizations. For MG--TH, the average ranges from $1.01$ to $1.06$ at $\Ra=10^3$ and $10^4$, and remains near $1.2$ at $\Ra=10^5$ for $N=64$ and $128$. The corresponding $N=256$ simulation did not reach the final time because of a linear solver convergence failure. For MG--AS, the averages remain between $1.04$ and $1.40$ across all nine configurations. At $\Ra=10^5$, the MG--AS iteration count increases only modestly, from $1.20$ at $N=64$ to $1.35$ at $N=256$. Overall, these results demonstrate weak mesh dependence of the outer Krylov iterations for the monolithic multigrid method.

The computational costs provide a more complete comparison of the preconditioners. At $N=256$ and $\Ra=10^4$, the average times per step are $2.138$, $1.770$, and $3.641$ seconds for AL--TH, MG--TH, and MG--AS, respectively, with corresponding outer Krylov iteration counts of $2.23$, $1.06$, and $1.40$. Thus, monolithic multigrid with Taylor--Hood elements reduces both iteration counts and computational time relative to AL--TH. With Alfeld--Sorokina elements, however, the larger algebraic system leads to higher computational costs despite the low iteration counts. At $N=128$ and $\Ra=10^5$, monolithic multigrid is faster with both discretizations, requiring $0.873$ seconds per step for MG--TH and $1.404$ seconds for MG--AS, compared with $1.606$ seconds for AL--TH.

We next investigate the influence of the Prandtl
number at $N=64$, with $\Delta t=1/64$ and $s=2$.
The selected values cover different regimes of
momentum and thermal diffusion: $\Pran=0.7$ and
$\Pran=7$ are representative of air and water near
room temperature, respectively, $\Pran=0.1$
represents low-Prandtl-number behavior associated with
thermally diffusive fluids, and $\Pran=1$ serves as
an intermediate reference case. Table~\ref{tab:parameter-results}
presents results for $\Ra=10^3$ and $10^4$ over
the range $0.1\leq\Pran\leq 7$. The monolithic multigrid method maintains average outer Krylov iteration counts between $1.00$ and $1.23$ for both finite element discretizations across all parameter combinations. In contrast, AL--TH exhibits greater sensitivity to the Prandtl number. At $\Ra=10^4$, its average iteration count increases from $3.52$ at $\Pran=0.1$ to $6.70$ at $\Pran=7$, whereas MG--TH and MG--AS maintain iteration counts close to one iteration per Newton step.

\begin{table}[!htbp]
\centering\small
\caption{Effect of the Prandtl number on solver performance
for $\Ra=10^3$ and $10^4$, with $N=64$, $s=2$,
$\Delta t=1/64$, and 16 MPI processes.
We report the average outer Krylov iterations per Newton
iteration ($\overline n_{\mathrm L}$) and the average
wall-clock time per step ($\overline t$, in seconds).}\label{tab:parameter-results}
\setlength{\tabcolsep}{5pt}
\begin{tabular}{rr rr rr rr}
\toprule
 & & \multicolumn{2}{c}{AL--TH} & \multicolumn{2}{c}{MG--TH} & \multicolumn{2}{c}{MG--AS}\\\cmidrule(lr){3-4}\cmidrule(lr){5-6}\cmidrule(lr){7-8}
$\Ra$ & $\Pran$ & $\overline n_{\mathrm L}$ & $\overline t$ & $\overline n_{\mathrm L}$ & $\overline t$ & $\overline n_{\mathrm L}$ & $\overline t$\\\midrule
$10^{3}$ & 0.1 & 1.63 & 0.186 & 1.00 & 0.371 & 1.03 & 0.303\\
$10^{3}$ & 0.7 & 1.55 & 0.169 & 1.02 & 0.159 & 1.05 & 0.303\\
$10^{3}$ & 1 & 1.71 & 0.155 & 1.02 & 0.148 & 1.04 & 0.296\\
$10^{3}$ & 7 & 5.41 & 0.227 & 1.02 & 0.229 & 1.05 & 0.426\\
$10^{4}$ & 0.1 & 3.52 & 0.405 & 1.01 & 0.295 & 1.01 & 0.481\\
$10^{4}$ & 0.7 & 3.00 & 0.307 & 1.09 & 0.205 & 1.10 & 0.330\\
$10^{4}$ & 1 & 3.21 & 0.292 & 1.10 & 0.223 & 1.11 & 0.331\\
$10^{4}$ & 7 & 6.70 & 0.347 & 1.15 & 0.210 & 1.23 & 0.345\\
\bottomrule\end{tabular}
\end{table}

The computational times do not necessarily follow the outer iteration counts, partly because of first-step assembly and compilation costs. For example, at $\Ra=10^3$ and $\Pran=0.1$, MG--TH requires $47.48$ seconds to complete the $128$ time steps, of which $33.59$ seconds are spent on the first step. This accounts for the relatively high average time per step reported in Table~\ref{tab:parameter-results}. Nevertheless, at $\Ra=10^4$, MG--TH is faster than AL--TH across all four Prandtl numbers. MG--AS also maintains low iteration counts at a generally higher computational cost.

\subsection{Stage dependence and computational cost}
\begin{figure}[!htbp]\centering
\includegraphics[width=\linewidth]{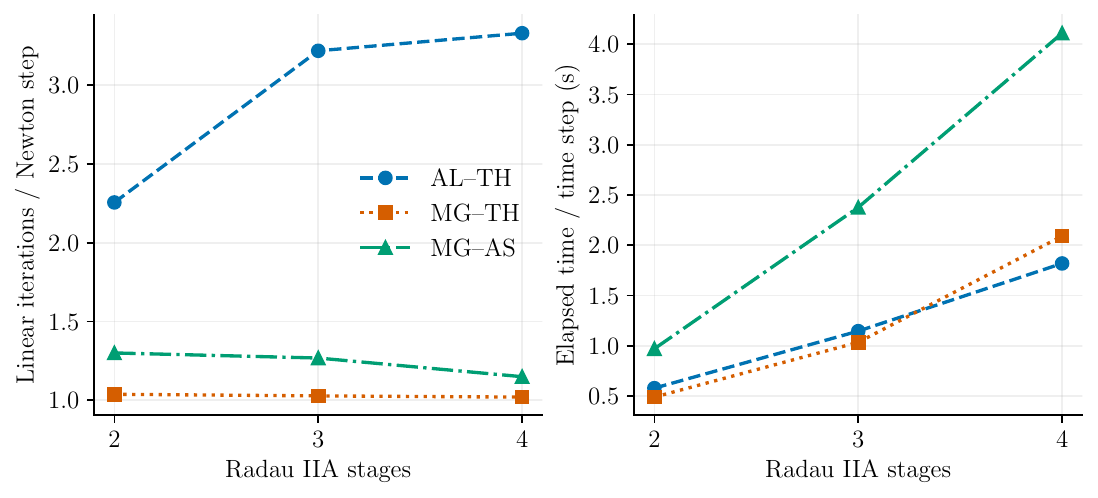}
\caption{Average outer Krylov iterations per Newton iteration
and run time per step as functions of the number
of Runge--Kutta stages for AL--TH, MG--TH, and MG--AS.
Results are obtained with $N=128$, $\Ra=10^4$, $\Pran=1$,
$\Delta t=1/128$, and 16 MPI processes.}\label{fig:stage-results}
\end{figure}

Table~\ref{tab:stage-results} presents the effect of increasing the number of Runge--Kutta stages at $\Ra=10^4$, $\Pran=1$, and $\Delta t=1/N$. Figure~\ref{fig:stage-results} illustrates the iteration counts and computational times at $N=128$. For AL--TH, increasing the number of stages from two to four raises the average outer Krylov iterations from $2.26$ to $3.33$ at $N=128$, and from $1.71$ to $2.60$ at $N=256$. In contrast, the monolithic multigrid method exhibits little variation in outer iteration counts with either finite element discretization. At $N=128$, MG--TH maintains averages of $1.04$, $1.03$, and $1.02$ for $s=2,3,4$, respectively. The corresponding averages for MG--AS decrease from $1.30$ to $1.15$. The average Newton iterations per time step for MG--AS also remain nearly constant, ranging from $0.89$ to $0.95$. At $N=256$, MG--TH maintains average outer Krylov
iteration counts between $1.01$ and $1.02$ across
all three stage counts. The corresponding range for MG--AS is $1.23$--$1.32$.

\begin{table}[!htbp]
\centering\small
\caption{Stage dependence with $\Ra=10^4$, $\Pran=1$, $\Delta t=1/N$, and 16 processes.}\label{tab:stage-results}
\setlength{\tabcolsep}{5pt}
\begin{tabular}{lrrrrr}
\toprule
Method & $N$ & $s$ & $\overline n_{\mathrm N}$ & $\overline n_{\mathrm L}$ & $\overline t$ (s)\\\midrule
AL--TH & 128 & 2 & 2.28 & 2.26 & 0.578\\
AL--TH & 128 & 3 & 2.27 & 3.22 & 1.146\\
AL--TH & 128 & 4 & 2.53 & 3.33 & 1.819\\
MG--TH & 128 & 2 & 0.84 & 1.04 & 0.489\\
MG--TH & 128 & 3 & 0.88 & 1.03 & 1.034\\
MG--TH & 128 & 4 & 0.82 & 1.02 & 2.092\\
MG--AS & 128 & 2 & 0.89 & 1.30 & 0.970\\
MG--AS & 128 & 3 & 0.91 & 1.27 & 2.377\\
MG--AS & 128 & 4 & 0.95 & 1.15 & 4.111\\
AL--TH & 256 & 2 & 2.10 & 1.71 & 1.662\\
AL--TH & 256 & 3 & 2.62 & 2.22 & 3.802\\
AL--TH & 256 & 4 & 2.54 & 2.60 & 5.760\\
MG--TH & 256 & 2 & 0.74 & 1.02 & 1.547\\
MG--TH & 256 & 3 & 0.77 & 1.01 & 3.194\\
MG--TH & 256 & 4 & 0.69 & 1.01 & 4.923\\
MG--AS & 256 & 2 & 0.81 & 1.32 & 3.281\\
MG--AS & 256 & 3 & 0.80 & 1.27 & 6.673\\
MG--AS & 256 & 4 & 0.86 & 1.23 & 12.150\\
\bottomrule\end{tabular}
\end{table}

Although the outer Krylov iteration counts remain low, computational costs increase with the number of Runge--Kutta stages. At $N=128$, the average step time for MG--AS increases from $0.970$ seconds for $s=2$ to $4.111$ seconds for $s=4$. The first step of the four-stage simulation alone requires $197.73$ seconds, reflecting substantial initial assembly and compilation costs. For $s=4$, MG--TH requires $2.092$ seconds per step, compared with $1.819$ seconds for AL--TH at $N=128$. At $N=256$, however, MG--TH is faster, requiring $4.923$ seconds per step compared with $5.760$ seconds for AL--TH. For MG--AS at $N=256$, increasing the stage count
from two to four raises the average step time from
$3.281$ to $12.150$ seconds, although the average
outer Krylov iterations per Newton iteration decrease
slightly from $1.32$ to $1.23$. These results illustrate how the relative computational costs vary with mesh resolution and the number of stages.

To further assess performance on larger problems, we consider two-stage simulations at $N=512$ with $\Delta t=1/512$. AL--TH and MG--TH require an average of $1.34$ and $1.04$ outer Krylov iterations per Newton iteration, with corresponding step times of $5.194$ and $5.403$ seconds, respectively. The coupled MG--TH system contains $5{,}255{,}176$ unknowns. These results demonstrate that both preconditioning approaches maintain low outer iteration counts for systems with more than five million unknowns, with comparable computational times.

\subsection{Parallel execution and reuse of solver components}
We investigate strong-scaling performance at $N=256$, $\Ra=10^4$, $\Pran=1$, $s=2$, and $\Delta t=1/256$, varying the number of MPI processes from four to sixteen. The parallel speedup and efficiency, relative to four processes, are defined by
\[
S_P=\frac{\overline t_4}{\overline t_P},
\qquad
E_P=\frac{S_P}{P/4}.
\]
As shown in Table~\ref{tab:parallel-results}, MG--AS achieves a speedup of $3.00$ and a parallel efficiency of $75.1\%$ on sixteen processes. The average outer Krylov iteration count remains approximately $1.33$ per Newton iteration. In comparison, AL--TH achieves a speedup of $1.92$, with average iteration counts ranging from $1.60$ to $1.71$. Thus, monolithic multigrid with Alfeld--Sorokina elements exhibits better strong-scaling performance over the tested process counts, although AL--TH retains shorter wall-clock times per step.

\begin{table}[!htbp]
\centering\small
\caption{Strong-scaling performance at $N=256$,
$\Ra=10^4$, $\Pran=1$, $s=2$, and $\Delta t=1/256$.
Speedup and parallel efficiency are computed relative
to the four-process run of each configuration.}\label{tab:parallel-results}
\setlength{\tabcolsep}{5pt}
\begin{tabular}{lrrrrr}
\toprule
Method & Processes & $\overline n_{\mathrm L}$ & $\overline t$ (s) & $S_P$ & $E_P$ (\%)\\\midrule
AL--TH & 4 & 1.60 & 3.186 & 1.00 & 100.0\\
AL--TH & 8 & 1.67 & 2.174 & 1.47 & 73.3\\
AL--TH & 16 & 1.71 & 1.662 & 1.92 & 47.9\\
MG--AS & 4 & 1.34 & 9.859 & 1.00 & 100.0\\
MG--AS & 8 & 1.34 & 5.622 & 1.75 & 87.7\\
MG--AS & 16 & 1.32 & 3.281 & 3.00 & 75.1\\
\bottomrule\end{tabular}
\end{table}

We next investigate the effect of lagging the Jacobian and preconditioner within the nonlinear solver. For $N=128$, $\Ra=10^5$, $\Pran=1$, $s=2$, and $\Delta t=1/128$, both the Jacobian and preconditioner are updated every five Newton iterations, with reuse allowed across time steps. Table~\ref{tab:reuse-results} compares these configurations with the corresponding baseline runs, where both are updated at every Newton iteration. Lagging reduces the average computational time per step by $32.1\%$ for MG--TH and $32.9\%$ for MG--AS, with only modest increases in Newton iterations and nearly unchanged outer Krylov iteration counts. These results demonstrate that Jacobian and preconditioner lagging can substantially reduce computational costs without significantly affecting linear solver convergence.

\begin{table}[!htbp]
\centering\small
\caption{Effect of Jacobian and preconditioner lagging on
solver performance at $N=128$, $\Ra=10^5$, $\Pran=1$,
$s=2$, $\Delta t=1/128$, and 16 MPI processes.
The update frequency applies to both the Jacobian
and preconditioner.}\label{tab:reuse-results}
\setlength{\tabcolsep}{5pt}
\begin{tabular}{llrrr}
\toprule
Method & Update frequency & $\overline n_{\mathrm N}$ & $\overline n_{\mathrm L}$ & $\overline t$ (s)\\\midrule
MG--TH & Every iteration & 1.41 & 1.24 & 0.873\\
MG--TH & Every 5 iterations & 1.52 & 1.27 & 0.593\\
MG--AS & Every iteration & 1.30 & 1.32 & 1.404\\
MG--AS & Every 5 iterations & 1.42 & 1.31 & 0.942\\
\bottomrule\end{tabular}
\end{table}

\subsection{Time resolution and cost}
We next investigate the tradeoff between temporal accuracy
and computational cost by varying the number of Runge--Kutta
stages and time steps. For the monolithic multigrid method
with both finite element discretizations, we fix $N=128$,
$\Ra=10^4$, $\Pran=1$, and $t_f=2$, using 16 MPI processes.
For each discretization, a four-stage Radau IIA solution
computed with 1024 time steps serves as the numerical
reference. An additional calculation with 512 time steps
is used to assess the accuracy of this reference.

We assess temporal accuracy by sampling the velocity and temperature fields on a fixed $33\times33$ grid at times $\mathcal T=\{1/8,1/4,1/2,1,2\}$. For either field $w$, we define the maximum sampled relative difference from the numerical reference solution as
\begin{equation}\label{eq:sampled-difference}
D_w=\max_{t\in\mathcal T}
\frac{
\left(\sum_j |w_h(x_j,t)-w_{\mathrm{ref}}(x_j,t)|^2\right)^{1/2}
}{
\left(\sum_j |w_{\mathrm{ref}}(x_j,t)|^2\right)^{1/2}
}.
\end{equation}
Here, $w_{\mathrm{ref}}$ denotes the corresponding reference solution computed with the same finite element discretization. Comparing the 512- and 1024-step reference calculations gives $(D_u,D_T)=(3.77\times10^{-7},1.22\times10^{-7})$ for Taylor--Hood and $(6.20\times10^{-7},1.89\times10^{-7})$ for Alfeld--Sorokina. These values provide an estimate of reference sensitivity when interpreting the temporal accuracy results.

\begin{table}[!htbp]
\centering\small
\caption{Temporal resolution and cost at $N=128$, $\Ra=10^4$, $\Pran=1$, $t_f=2$, and 16 processes. $\Delta t=2/n_t$. The sampled relative differences $D_u$ and $D_T$ are defined in \eqref{eq:sampled-difference}.}\label{tab:timestep-results}
\setlength{\tabcolsep}{5pt}
\begin{tabular}{lrrrrr}
\toprule
Pair & $s$ & $n_t$ & Time (s) & $D_u$ & $D_T$\\\midrule
TH & 2 & 128 & 80.74 & $2.222\times10^{-3}$ & $1.126\times10^{-3}$\\
TH & 2 & 64 & 45.58 & $1.253\times10^{-2}$ & $6.794\times10^{-3}$\\
TH & 3 & 128 & 157.58 & $5.405\times10^{-5}$ & $1.788\times10^{-5}$\\
TH & 3 & 64 & 94.96 & $1.451\times10^{-3}$ & $4.916\times10^{-4}$\\
TH & 4 & 128 & 265.13 & $4.351\times10^{-6}$ & $1.399\times10^{-6}$\\
TH & 4 & 64 & 157.96 & $5.903\times10^{-5}$ & $2.034\times10^{-5}$\\
AS & 2 & 128 & 148.61 & $2.222\times10^{-3}$ & $1.126\times10^{-3}$\\
AS & 2 & 64 & 91.14 & $1.253\times10^{-2}$ & $6.794\times10^{-3}$\\
AS & 3 & 128 & 319.85 & $5.419\times10^{-5}$ & $1.793\times10^{-5}$\\
AS & 3 & 64 & 186.48 & $1.451\times10^{-3}$ & $4.916\times10^{-4}$\\
AS & 4 & 128 & 529.54 & $1.245\times10^{-6}$ & $3.778\times10^{-7}$\\
AS & 4 & 64 & 322.50 & $5.908\times10^{-5}$ & $2.033\times10^{-5}$\\
\bottomrule\end{tabular}
\end{table}

\begin{figure}[!htbp]\centering
\includegraphics[width=\linewidth]{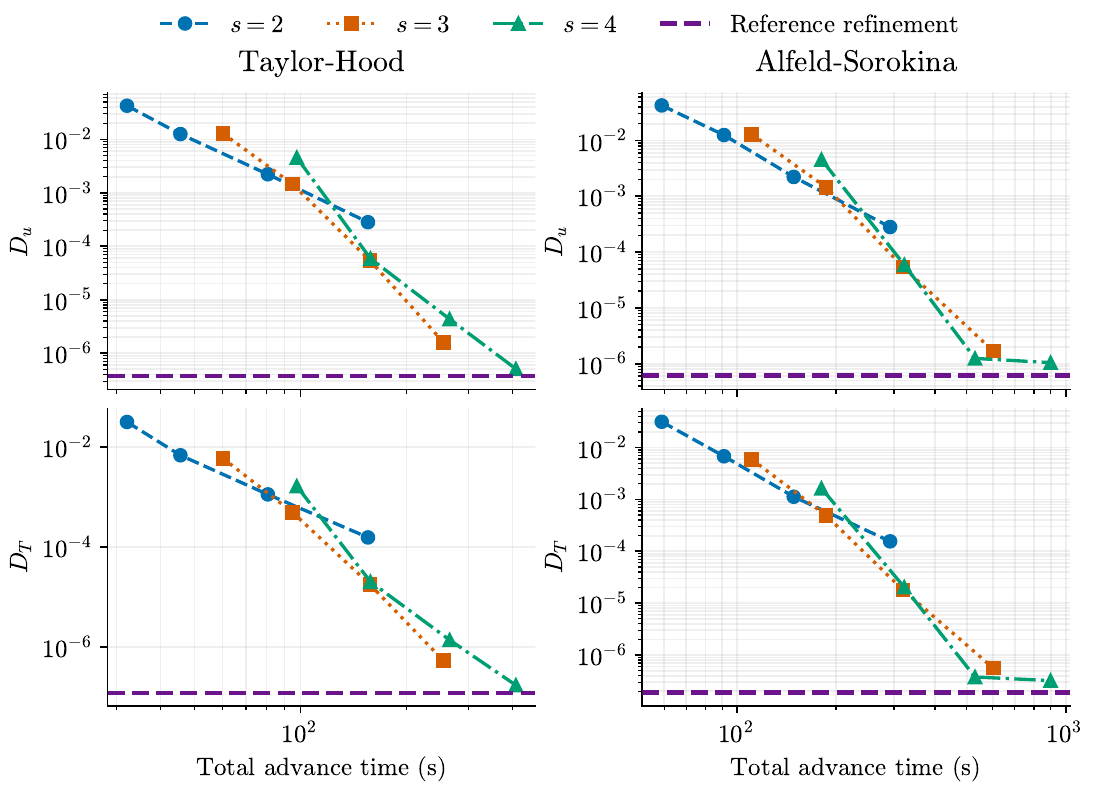}
\caption{Temporal accuracy versus computational cost for Taylor--Hood (left) and Alfeld--Sorokina (right) discretizations. The top and bottom panels show the
sampled relative differences in velocity ($D_u$) and
temperature ($D_T$), respectively, as functions of total
computational time. Results are shown for two-, three-,
and four-stage Radau IIA methods with 32, 64, 128,
and 256 time steps. Horizontal dashed lines
indicate the differences between the 512- and 1024-step
reference solutions.}\label{fig:timestep-results}
\end{figure}

Table~\ref{tab:timestep-results} and Figure~\ref{fig:timestep-results} demonstrate the advantages of higher-order time integration. For Taylor--Hood elements, the three-stage method with 128 time steps achieves velocity and temperature differences approximately $5.2$ and $8.7$ times smaller, respectively, than the two-stage method with 256 steps, at nearly identical computational costs ($157.58$ versus $155.44$ seconds). For Alfeld--Sorokina elements, three stages with 128 steps and four stages with 64 steps achieve comparable accuracy, requiring $319.85$ and $322.50$ seconds, respectively. Both configurations provide smaller temporal differences than the two-stage method with 256 steps, which requires $292.05$ seconds. These results show that higher-stage Radau IIA methods can achieve greater temporal accuracy at comparable computational cost.

\subsection{Comparison with incremental pressure correction}
\label{subsec:radau-pressure-correction}

We next compare the fully implicit Radau IIA methods with a
BDF2-based incremental pressure correction scheme
(BDF2--IPCS),
a widely used projection method for incompressible flows
\cite{GuermondMinevShen2006}. Unlike fully implicit Runge--Kutta
methods, which solve the coupled velocity--pressure--temperature
system, projection methods decouple the velocity and pressure
calculations. This reduces the computational cost of
individual time steps, but higher-order implicit
integration may achieve greater accuracy with
substantially fewer steps.

The BDF2--IPCS implementation follows the nonrotational
incremental pressure-correction formulation, with a
backward-Euler startup step. The tentative velocity uses
extrapolated advection velocity and temperature,
$\mathbf w^n=2\mathbf u^n-\mathbf u^{n-1}$ and
$T^{*,n}=2T^n-T^{n-1}$. Defining the pressure increment
$\phi^{n+1}=p^{n+1}-p^n$ and $\mathbf g=(0,-1)^T$,
the velocity--pressure update is given by
\begin{align}
\frac{3\widetilde{\mathbf u}^{n+1}
-4\mathbf u^n+\mathbf u^{n-1}}{2\Delta t}
+(\mathbf w^n\cdot\nabla)\widetilde{\mathbf u}^{n+1}
-\Pran\Delta\widetilde{\mathbf u}^{n+1}
&\nonumber\\
+\Pran\nabla p^n
+\Ra\Pran T^{*,n}\mathbf g
&=\mathbf 0,
\label{eq:ipcs-rbc-predictor}\\
(\nabla\phi^{n+1},\nabla q)
&=-\frac{3}{2\Pran\Delta t}
(\nabla\cdot\widetilde{\mathbf u}^{n+1},q),
\label{eq:ipcs-rbc-pressure}\\
(\mathbf u^{n+1},\mathbf v)
&=(\widetilde{\mathbf u}^{n+1},\mathbf v)
-\frac{2\Pran\Delta t}{3}
(\nabla\phi^{n+1},\mathbf v).
\label{eq:ipcs-rbc-correction}
\end{align}
The temperature is subsequently advanced using an implicit
BDF2 advection--diffusion solve with the corrected velocity.
The pressure-increment problem is solved with homogeneous
Neumann conditions and its constant nullspace removed.

The comparison uses the two-dimensional benchmark with
Taylor--Hood elements, $N=128$, $\Ra=10^4$, $\Pran=1$,
and $t_f=2$, using 16 MPI processes.
Both methods use the same initial conditions and boundary
conditions. We measure the sampled relative differences
$D_u$ and $D_T$ defined in \eqref{eq:sampled-difference},
using the same four-stage, 1024-step Radau IIA numerical
reference. The tentative-velocity and temperature systems are solved
using FGMRES, with geometric multigrid preconditioning.
On each noncoarse level, five GMRES iterations are applied
with block-Jacobi preconditioning and local ILU factorizations.
The pressure-increment Poisson problem is solved using CG
preconditioned by geometric multigrid, with three Chebyshev
smoothing iterations using Jacobi preconditioning on each
noncoarse level. All three systems use direct MUMPS coarse
solves, relative and absolute tolerances of $10^{-8}$ and
$10^{-10}$, respectively, and a maximum of 200 Krylov
iterations. The velocity correction is computed through a
finite element mass-matrix solve using CG with Jacobi
preconditioning and relative and absolute tolerances
of $10^{-10}$ and $10^{-12}$, respectively.
Radau IIA uses the monolithic multigrid solver
described above.

\begin{table}[htbp]
\centering
\caption{Temporal accuracy and computational cost of
BDF2--IPCS and fully implicit Radau IIA at $N=128$,
$\Ra=10^4$, $\Pran=1$, and $t_f=2$.
Total wall-clock time includes time advancement and
first-step setup costs. The sampled relative differences
are defined in \eqref{eq:sampled-difference}.}
\label{tab:radau-ipcs-work}
\begingroup
\small
\setlength{\tabcolsep}{3pt}
\begin{tabular}{lrrrrrr}
\toprule
Method & Steps & $\Delta t$ & Wall/step & Total wall & $D_u$ & $D_T$ \\
 & & & (s) & (s) & & \\
\midrule
BDF2--IPCS & 1024 & $1/512$ & 0.06240 & 63.90 & $9.096\times10^{-3}$ & $6.751\times10^{-3}$ \\
BDF2--IPCS & 2048 & $1/1024$ & 0.05220 & 106.92 & $2.310\times10^{-3}$ & $1.719\times10^{-3}$ \\
BDF2--IPCS & 4096 & $1/2048$ & 0.04724 & 193.48 & $5.865\times10^{-4}$ & $4.368\times10^{-4}$ \\
Radau IIA, $s=2$ & 256 & $1/128$ & 0.51075 & 130.75 & $2.820\times10^{-4}$ & $1.556\times10^{-4}$ \\
Radau IIA, $s=3$ & 128 & $1/64$ & 1.21443 & 155.45 & $5.405\times10^{-5}$ & $1.788\times10^{-5}$ \\
Radau IIA, $s=4$ & 64 & $1/32$ & 2.46745 & 157.92 & $5.903\times10^{-5}$ & $2.034\times10^{-5}$ \\
\bottomrule
\end{tabular}
\endgroup
\end{table}

The BDF2--IPCS results in
Table~\ref{tab:radau-ipcs-work} exhibit approximately
second-order decay of the sampled reference differences
under time-step refinement. Using
\[
p_w=\log_2\left(
\frac{D_w(\Delta t)}{D_w(\Delta t/2)}
\right),
\]
the observed orders for successive refinements are
$1.98$ and $1.98$ for velocity, and $1.97$ and $1.98$
for temperature. These results provide numerical evidence
of second-order temporal convergence in the sampled
quantities for this benchmark, although they do not
establish a general convergence result for the coupled
pressure-correction scheme.

Table~\ref{tab:radau-ipcs-work} demonstrates the potential
advantage of fully implicit Runge--Kutta integration.
Compared with BDF2--IPCS using 4096 time steps,
the two-, three-, and four-stage Radau IIA methods
require only 256, 128, and 64 steps, respectively,
and achieve smaller sampled differences in both
velocity and temperature. Their total computational
times are also reduced by approximately $18$--$32\%$.
In particular, the three-stage method reduces the
velocity and temperature differences by factors of
approximately $10.9$ and $24.4$, respectively, and
requires $19.7\%$ less computational time.

Although each fully implicit Runge--Kutta step is more
expensive, its higher temporal accuracy can substantially
reduce the total computational cost. These results
demonstrate that fully implicit integration can outperform
a projection method when sufficiently accurate transient
solutions are required.

\subsection{Three-dimensional experiments}
\label{sec:three_dimensional_rbc}

We extend our numerical study to three-dimensional
Rayleigh--B\'enard convection, examining the performance
of the monolithic multigrid method, the effect of
coarse-grid solvers, and parallel scalability.
We consider both Taylor--Hood and enriched
Guzm\'an--Neilan finite element discretizations.
All experiments are performed on the Magnolia cluster
at the Mississippi Center for Supercomputing Research
(MCSR).

\paragraph{Computational setting}
We consider the unit cube $\Omega=(0,1)^3$, with $y$ as
the vertical coordinate. No-slip velocity conditions and
the thermal boundary conditions~\eqref{eq:boundary-conditions}
are imposed, with the initial temperature given by
\eqref{eq:initial-conditions}, independent of $z$.
The fluid is initially at rest. Unless otherwise specified,
we take $\Ra=10^4$ and $\Pran=1$.
We use Taylor--Hood elements (continuous quadratic velocity
and linear pressure) with continuous piecewise linear
temperature. Time integration employs the two-stage
Radau IIA method with $\Delta t=1/N$, where $N$ denotes
the number of mesh subdivisions in each coordinate direction.
All Runge--Kutta stages and physical fields are solved
simultaneously.

The coupled nonlinear systems are solved using an inexact Newton method with FGMRES and monolithic geometric multigrid preconditioning. The multigrid method employs vertex-centered overlapping Vanka relaxation, with local patch systems solved by sparse LU factorization. For the coarse-grid solve, we compare direct MUMPS factorization with PETSc \cite{petsc-user-ref} Telescope, which reduces the MPI communicator before applying MUMPS. We report total Newton and outer Krylov iterations, together with the average computational time per step. For the $N=64$ strong-scaling study, we also evaluate
the speedup excluding the first time step to assess
the influence of initial setup costs.

\begin{figure}[!htbp]\centering
\begin{minipage}[b]{.53\linewidth}\centering
\includegraphics[width=\linewidth,height=5.0cm,keepaspectratio]{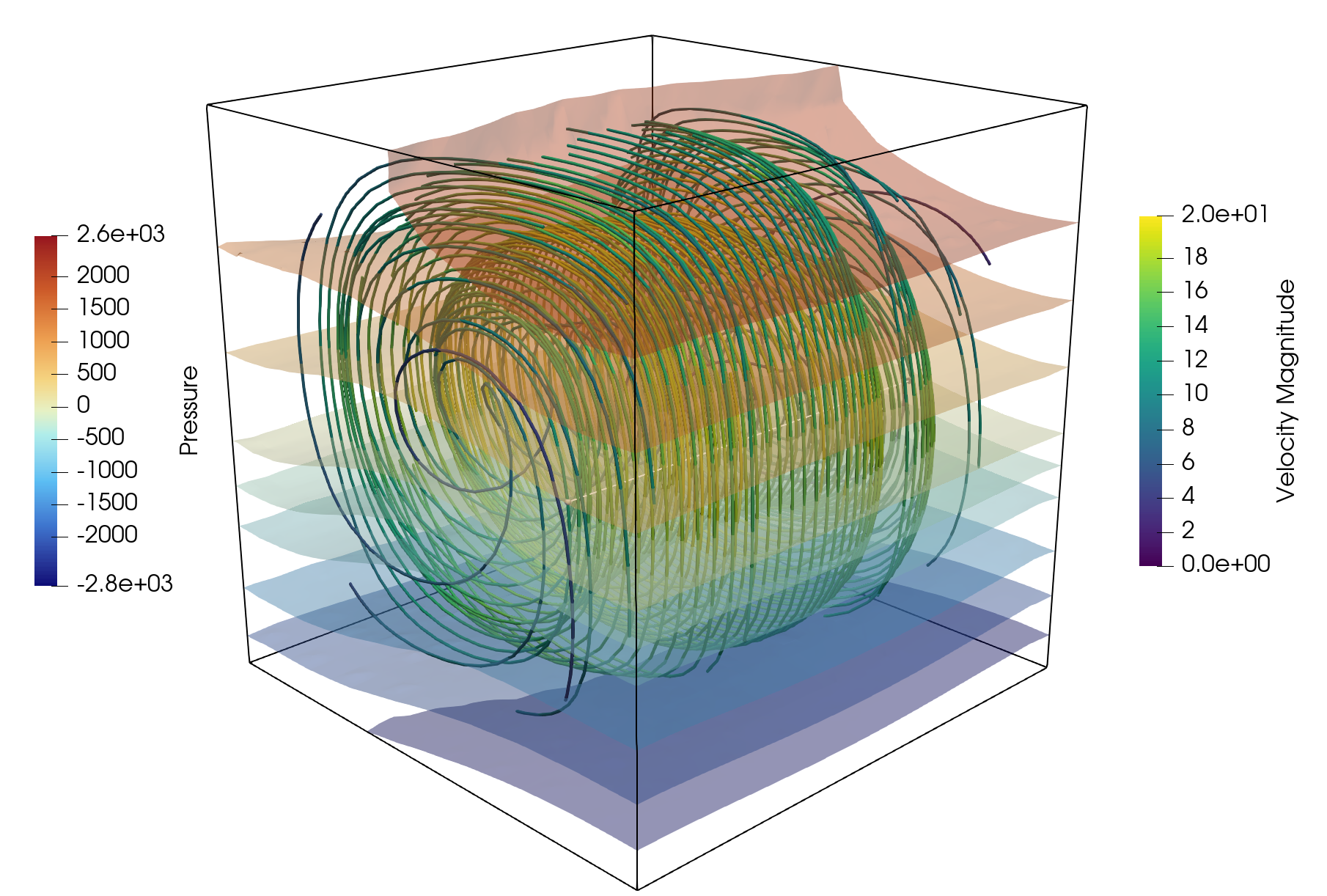}\\
\small (a) Velocity streamlines and pressure slices
\end{minipage}\hfill
\begin{minipage}[b]{.45\linewidth}\centering
\includegraphics[width=\linewidth,height=5.0cm,keepaspectratio]{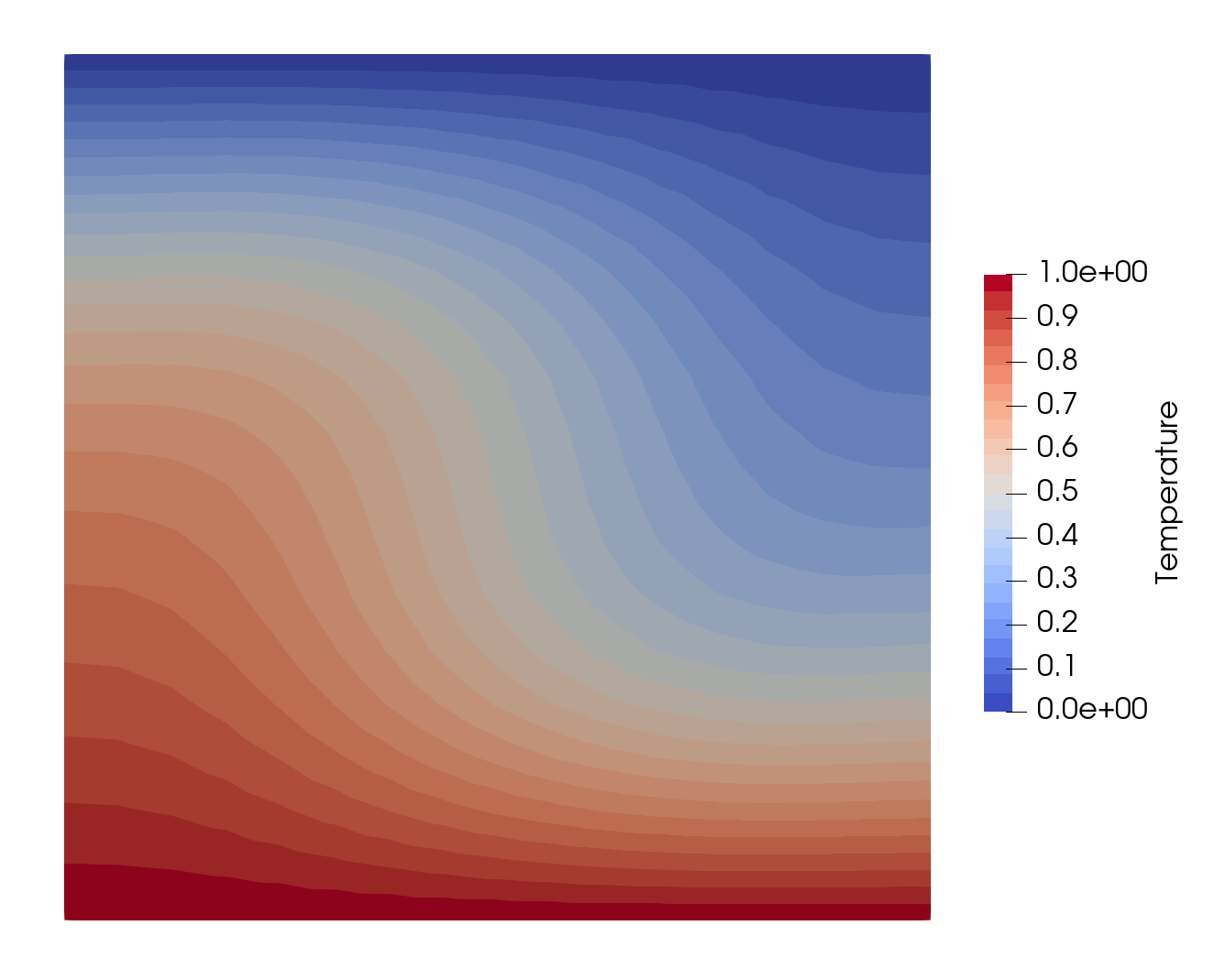}\\
\small (b) Temperature slice
\end{minipage}
\caption{Three-dimensional Rayleigh--B\'enard convection
with Taylor--Hood elements at $N=16$, $\Ra=10^4$,
$\Pran=1$, and $t=2$, using two-stage Radau IIA
with $\Delta t=1/16$.
(a) Velocity streamlines colored by velocity magnitude,
together with horizontal pressure slices.
(b) Temperature distribution on a cross-sectional slice.}
\label{fig:three-d-snapshots}

\end{figure}
Figure~\ref{fig:three-d-snapshots} illustrates the three-dimensional circulating flow and the associated temperature distribution. We next examine the velocity divergence and wall heat fluxes to assess mass conservation and thermal transport.

\paragraph{Physical observables and discrete incompressibility}
We now assess the three-dimensional solutions using the
$L^2$ norms of velocity, temperature, and velocity
divergence, together with the Nusselt numbers at the
horizontal boundaries,
\begin{equation}
Nu_{\mathrm{bot}}
=-\int_{y=0}\partial_y T_h\,\mathrm{d}S,
\qquad
Nu_{\mathrm{top}}
=-\int_{y=1}\partial_y T_h\,\mathrm{d}S.
\label{eq:three_d_nusselt}
\end{equation}

The Taylor--Hood simulation shown in
Figure~\ref{fig:three-d-snapshots} reaches $t_f=2$
after 32 time steps, yielding
$\|u_h\|_{L^2}=11.43184$ and
$\|\nabla\cdot u_h\|_{L^2}=1.73991$.
The corresponding Nusselt numbers are
$Nu_{\mathrm{bot}}=2.00180$ and
$Nu_{\mathrm{top}}=2.00268$, indicating approximately
twice the conductive heat flux and close agreement
between the two boundaries.

To assess the effect of spatial discretization on mass
conservation, we compare Taylor--Hood elements with
the enriched Guzm\'an--Neilan velocity space
(\texttt{GNH1div}) of degree two, paired with
continuous piecewise linear pressure on the
Alfeld split~\cite{GuzmanNeilan2018,BrubeckKirby2025}.
The temperature is approximated using continuous
piecewise linear elements. The two discretizations
are compared after one time step, at $t=1/16$.

\begin{table}[t]
\centering
\caption{Three-dimensional solution diagnostics at $N=16$, $t=1/16$,
$\Ra=10^4$, $\Pran=1$, with four MPI ranks and two-stage Radau IIA.  All
runs use $\Delta t=1/16$.  Here TH denotes Taylor--Hood and GN denotes
the enriched \texttt{GNH1div}--Alfeld pair.}
\label{tab:three_d_discrete_diagnostics}
\small
\begin{tabular}{lrrrrr}
\toprule
Discretization & $\|\mathbf u_h\|_{L^2}$ & $\|T_h\|_{L^2}$ &
$\|\nabla\cdot\mathbf u_h\|_{L^2}$ & $Nu_{\mathrm{bot}}$ & $Nu_{\mathrm{top}}$ \\
\midrule
TH, direct    & 4.909641 & 0.575976 & $7.8415\times10^{-1}$ & 1.117931 & 1.117912 \\
GN, direct    & 4.898039 & 0.575980 & $2.2931\times10^{-8}$ & 1.117475 & 1.117452 \\
GN, Tel.~$\times2$ & 4.898039 & 0.575980 & $2.2911\times10^{-8}$ & 1.117475 & 1.117452 \\
GN, Tel.~$\times4$ & 4.898039 & 0.575980 & $2.2911\times10^{-8}$ & 1.117475 & 1.117452 \\
\bottomrule
\end{tabular}
\end{table}

Table~\ref{tab:three_d_discrete_diagnostics} shows close agreement between the two discretizations in the velocity and temperature norms and the wall Nusselt numbers. However, the velocity divergence norm decreases from $7.84\times10^{-1}$ for Taylor--Hood to approximately $2.3\times10^{-8}$ for the Guzm\'an--Neilan discretization. This substantial reduction reflects the divergence-free properties of the enriched velocity--pressure pair and its ability to achieve nearly exact discrete mass conservation.

\paragraph{Coarse-grid treatment and parameter dependence}

On the $N=16$ Taylor--Hood mesh, direct MUMPS
and PETSc Telescope require 58 Newton and 72
outer Krylov iterations over 32 time steps,
with comparable computational times.

We also investigate the dependence of solver performance on the Rayleigh and Prandtl numbers. Table~\ref{tab:three_d_parameters} presents six eight-step simulations using Taylor--Hood elements at $N=16$, with $\Ra=10^3,10^4$ and $\Pran=0.7,1,7$. At $\Ra=10^3$, the solver requires between 16 and 19 Newton iterations and 20 to 27 outer Krylov iterations across the three Prandtl numbers. At $\Ra=10^4$, increasing $\Pran$ from $0.7$ to $7$ raises the Newton iterations from 27 to 62 and the outer Krylov iterations from 41 to 98, with average step times increasing from $53.06$ to $123.68$ seconds. All six simulations complete successfully using the same monolithic multigrid preconditioner.

\begin{table}[t]
\centering
\caption{Parameter dependence of the three-dimensional
Taylor--Hood solver at $N=16$, with four MPI processes,
$\Delta t=1/16$, and eight time steps.
The coarse-grid problem is solved using MUMPS.}
\label{tab:three_d_parameters}
\small
\begin{tabular}{rrrrrr}
\toprule
$\Ra$ & $\Pran$ & $n_{\mathrm N}$ & $n_{\mathrm L}$
& $\overline n_{\mathrm L}$ & Time/step (s) \\
\midrule
$10^3$ & 0.7 & 18 & 20 & 1.11 & 31.48 \\
$10^3$ & 1.0 & 19 & 21 & 1.11 & 36.37 \\
$10^3$ & 7.0 & 16 & 27 & 1.69 & 32.69 \\
\midrule
$10^4$ & 0.7 & 27 & 41 & 1.52 & 53.06 \\
$10^4$ & 1.0 & 28 & 42 & 1.50 & 60.11 \\
$10^4$ & 7.0 & 62 & 98 & 1.58 & 123.68 \\
\bottomrule
\end{tabular}
\end{table}

\paragraph{Parallel performance}

We investigate the parallel performance of the monolithic
multigrid method using three-dimensional Taylor--Hood
discretizations. Table~\ref{tab:three_d_parallel} presents
results for $N=16,32,64$ at different MPI process counts.
These experiments examine solver performance as the problem
size increases. Comparisons at fixed mesh resolution provide strong-scaling results.

\begin{table}[t]
\centering
\caption{Three-dimensional Taylor--Hood scale-up experiments
for $\Ra=10^4$, $\Pran=1$, and two-stage Radau IIA.
Each simulation advances four time steps with $\Delta t=1/N$.
The $N=32$ runs use the mesh hierarchy $8\to16\to32$.
Tel.\ denotes a Telescope coarse-grid solve,
and $P_c$ is the number of MPI ranks used for coarse-grid
factorization. The reported iteration counts are totals
over the four time steps.}
\label{tab:three_d_parallel}
\small
\setlength{\tabcolsep}{4.5pt}
\begin{tabular}{rrlrrrr}
\toprule
$N$ & MPI ranks & Coarse solver & $P_c$
& $n_{\mathrm N}$ & $n_{\mathrm L}$
& Time/step (s) \\
\midrule
16 & 4   & Direct           & 4 & 16 & 25 &  70.56 \\
16 & 4   & Tel.~$\times4$  & 1 & 16 & 25 &  78.06 \\
\midrule
32 & 32  & Tel.~$\times4$  & 8 & 14 & 23 &  90.73 \\
32 & 64  & Tel.~$\times8$  & 8 & 14 & 23 &  58.07 \\
32 & 128 & Tel.~$\times16$ & 8 & 14 & 23 &  30.71 \\
\midrule
64 & 64  & Tel.~$\times8$  & 8 & 12 & 18 & 351.02 \\
64 & 128 & Tel.~$\times16$ & 8 & 12 & 18 & 173.52 \\
64 & 256 & Tel.~$\times32$ & 8 & 12 & 18 &  98.31 \\
\bottomrule
\end{tabular}
\end{table}
For $N=32$, increasing the MPI process count from 32 to 128
reduces the average time per step from 90.73 to 30.71 seconds,
giving a speedup of 2.95 and approximately 74\% parallel
efficiency. All three configurations require 14 Newton
and 23 outer Krylov iterations.
Stronger scaling is observed for $N=64$ using 64, 128,
and 256 MPI processes. Both scaling series maintain
eight MPI processes for the coarse-grid factorization. As shown in Table~\ref{tab:three_d_parallel}, the average
time per step decreases from $351.02$ to $173.52$ and
$98.31$ seconds, respectively. All three configurations
require 12 Newton and 18 outer Krylov iterations.
Relative to the 64-process configuration, the speedup
and parallel efficiency on 256 processes are
\begin{equation}
S_{64\rightarrow256}
=\frac{351.02}{98.31}\approx3.57,
\qquad
E_{64\rightarrow256}
=\frac{S_{64\rightarrow256}}{4}\approx0.893.
\label{eq:three_d_strong_scaling}
\end{equation}
Excluding the first time step yields a comparable speedup
of $3.41$. The velocity and temperature norms, velocity
divergence, and Nusselt numbers agree across all three
configurations to the reported precision. These results
demonstrate favorable strong-scaling performance for
the fully coupled three-dimensional systems, achieving
approximately $89\%$ parallel efficiency.

\paragraph{Enriched elements}

The enriched Guzm\'an--Neilan (\texttt{GNH1div}--Alfeld)
discretization is more computationally demanding than
Taylor--Hood at the same nominal mesh resolution.
At $N=16$, the coupled two-stage system contains
$395\,628$ unknowns, with approximately $88.7$ million
nonzero entries in the assembled Jacobian.
For $\Ra=10^4$ and $\Pran=1$, the first time step
converges in four Newton and eleven outer Krylov
iterations for all three coarse-grid configurations.
The corresponding computational times are $286.28$
seconds with direct MUMPS and $294.04$ and $300.35$
seconds with Telescope reduction factors of two and
four, respectively. The solution diagnostics agree
across these configurations, as shown in
Table~\ref{tab:three_d_discrete_diagnostics}.

An eight-step simulation at $\Ra=10^3$ and $\Pran=1$
also reaches $t=0.5$, requiring 15 Newton and 169
outer Krylov iterations. The velocity norm is
$3.27\times10^{-4}$, the divergence norm is
$1.02\times10^{-6}$, and both Nusselt numbers are
approximately $1.00001$, consistent with a state
close to thermal conduction. For the more strongly
driven $\Ra=10^4$ problem, the enriched configurations
encounter inner linear-solver convergence difficulties
beyond the first time step. These experiments establish
the feasibility of solving fully coupled Runge--Kutta
systems with enriched Guzm\'an--Neilan elements.
Improving multistep robustness through Vanka patch
construction, Krylov restarting, and coarse-grid
stabilization remains a direction for future work.

\section{Conclusions and future work}
\label{sec:conclusions}

In this paper, we have investigated augmented Lagrangian block
preconditioning and monolithic geometric multigrid
for fully implicit Radau IIA discretizations of
Rayleigh--B\'enard convection. The numerical experiments
demonstrate the effectiveness of both approaches for
solving the coupled velocity--pressure--temperature
systems. In particular, monolithic multigrid with
vertex-centered Vanka relaxation maintains low outer
Krylov iteration counts across mesh refinements,
physical parameter variations, and increasing numbers
of Runge--Kutta stages.

The temporal accuracy studies demonstrate an important
advantage of higher-order implicit integration.
Although individual Radau IIA steps are more expensive,
higher-stage methods can achieve greater accuracy with
substantially fewer time steps. At the $10^{-3}$ sampled-difference threshold, the tested
Radau IIA configurations achieve smaller velocity and
temperature differences than BDF2--IPCS and reduce total
computational time by 18--32\%.
These results highlight the importance of considering
temporal accuracy and total computational cost together.

The three-dimensional experiments further demonstrate
the effectiveness of monolithic multigrid on parallel
systems. For the $64^3$ Taylor--Hood discretization,
increasing the number of MPI processes from 64 to 256
achieves a speedup of $3.57$, corresponding to $89.3\%$
parallel efficiency, with unchanged nonlinear and
Krylov iteration counts. The enriched Guzm\'an--Neilan
discretization also achieves nearly divergence-free
velocity fields, demonstrating the benefits of
compatible finite element spaces for discrete
mass conservation.

Future work will focus on improving the efficiency
of Vanka relaxation through patch coloring, batched
factorizations, and GPU acceleration. Further investigation of efficient patch solvers
and coarse-grid stabilization is needed to improve
the robustness of enriched discretizations over
longer simulations. Extending
the accuracy and performance studies to higher
Rayleigh numbers and larger three-dimensional
problems is another important direction of our future work.

\section*{Acknowledgments}
The authors acknowledge the Mississippi Center for Supercomputing Research for access to the Magnolia cluster.
\bibliographystyle{siamplain}
\bibliography{references}
\end{document}

%% file: ex_shared.tex
\usepackage{amsfonts}
\usepackage{graphicx}
\usepackage{epstopdf}
\usepackage{algorithmic}
\ifpdf
  \DeclareGraphicsExtensions{.eps,.pdf,.png,.jpg}
\else
  \DeclareGraphicsExtensions{.eps}
\fi

\newsiamremark{remark}{Remark}
\newsiamremark{hypothesis}{Hypothesis}
\crefname{hypothesis}{Hypothesis}{Hypotheses}
\newsiamthm{claim}{Claim}
\newsiamremark{fact}{Fact}
\crefname{fact}{Fact}{Facts}

\headers{Time-Dependent Rayleigh--B\'enard Convection}{Ahsan Ali and Robert C. Kirby}

\title{Scalable Fully Implicit Runge--Kutta Solvers for
Time-Dependent Rayleigh--B\'enard Convection\thanks{
\funding{This work was supported by National Science Foundation grant OAC-2410408.}}}

\author{
Ahsan Ali\thanks{
School of Mathematics and Natural Sciences,
University of Southern Mississippi,
Hattiesburg, MS 39406, USA
(\email{ahsan.ali@usm.edu}).
}
\and
Robert C. Kirby\thanks{
Department of Mathematics,
Baylor University,
Waco, TX 76706, USA
(\email{robert\_kirby@baylor.edu}).
}
}

\usepackage{amsopn}
